\documentclass[preprint,12pt]{elsarticle}

\usepackage{amssymb}

\usepackage[utf8]{inputenc}
\usepackage{lmodern} 
\usepackage[T1]{fontenc}
\usepackage{microtype}

\usepackage{graphicx}
\usepackage{float}
\usepackage[hypcap]{caption}
\usepackage{subcaption}
\usepackage{algorithm2e}
\usepackage{array}
\usepackage{algorithmic}

\usepackage{textcomp}
\usepackage{booktabs}
\usepackage{array}

\usepackage{amsfonts}

\usepackage{geometry}
\usepackage{tikz}
\usetikzlibrary{trees,positioning,shapes}

\usepackage{mathtools}
\usepackage{amssymb}
\usepackage{amsthm}
\usepackage{amsmath}
\usepackage{xcolor, colortbl}
\usepackage{appendix}

\usepackage[english]{babel}
\usepackage{url}

\usepackage{listings}

\lstdefinestyle{custompython}{
  belowcaptionskip=1\baselineskip,
  breaklines=true,
  frame=L,
  xleftmargin=\parindent,
  language=python,
  showstringspaces=false,
  basicstyle=\footnotesize\ttfamily,
  keywordstyle=\bfseries\color{green!40!black},
  commentstyle=\itshape\color{purple!40!black},
  identifierstyle=\color{blue},
  stringstyle=\color{orange},
}
\journal{Applied Soft Computing}

\begin{document}

\begin{frontmatter}



\title{Graph grammars and Physics Informed Neural Networks for simulating of pollution propagation on Spitzbergen}


\author[inst2]{Maciej Sikora}

\author[inst1]{Albert Oliver-Serra}


\author[inst2]{Leszek Siwik}

\author[inst3]{\\ Natalia Leszczy\'nska}

\author[inst4]{Tomasz Maciej Ciesielski}

\author[inst5,inst6,inst7]{\\ Eirik Valseth}

\author[inst2]{Jacek Leszczy\'nski}

\author[inst8]{Anna Paszy\'nska}

\author[inst2]{Maciej Paszy\'nski}

\affiliation[inst2]{organization={AGH University of Krakow, Faculty of Computer Science},
            addressline={Al. Mickiewicza 30}, 
            city={Krak\'ow},
            postcode={30-059}, 
            country={Poland}}

\affiliation[inst1]{organization={University Institute of Intelligent Systems and Numeric Applications in Engineering, University of Las Palmas de Gran Canaria (ULPGC)},
            city={Las Palmas de Gran Canaria},
            country={Spain}}

\affiliation[inst3]{organization={Medical University of Silesia-Katowice, Faculty of Medical Sciences}, 
city={Katowice},
country={Poland}}

\affiliation[inst4]{organization={The University Centre in Svalbard},
            city={Longyearbyen},
            country={Norway}}

\affiliation[inst5]{organization={Norwegian University of Life Sciences},
            city={Ås},
            country={Norway}}

\affiliation[inst6]{organization={The Oden Institute for Computational Engineering Sciences, , The University of Texas},
            city={Austin, Texas},
            country={USA}}

\affiliation[inst7]{organization={Simula Research Laboratory},
            city={Oslo},
            country={Norway}}

\affiliation[inst8]{organization={Jagiellonian University, Faculty of Physics, Astronomy and Applied Computer Science},
            addressline={\L{}ojasiewicza 11}, 
            city={Krak\'ow},
            postcode={30-348}, 
            country={Poland}}

\begin{abstract}

In this paper, we present two computational methods for performing simulations of pollution propagation described by advection-diffusion equations. The first method employs graph grammars to describe the generation process of the computational mesh used in simulations with the meshless solver of the three-dimensional finite element method. 
The graph transformation rules express the three-dimensional Rivara longest-edge refinement algorithm.
This solver is used for an exemplary application: performing three-dimensional simulations of pollution generation by the coal-burning power plant and its propagation in the city of Longyearbyen, the capital of Spitsbergen.
The second computational code is based on the Physics Informed Neural Networks method. It is used to calculate the dissipation of the pollution along the valley in which the city of Longyearbyen is located. We discuss the instantiation and execution of the PINN method using Google Colab implementation.
We discuss the benefits and limitations of the PINN implementation.
\end{abstract}

\begin{keyword}
Graph grammar model \sep Physics Informed Neural Networks \sep Pollution propagation simulations 
\end{keyword}

\end{frontmatter}

\section{Introduction}

 \begin{figure}[ht]
    \centering
    \includegraphics[width=0.8\linewidth]{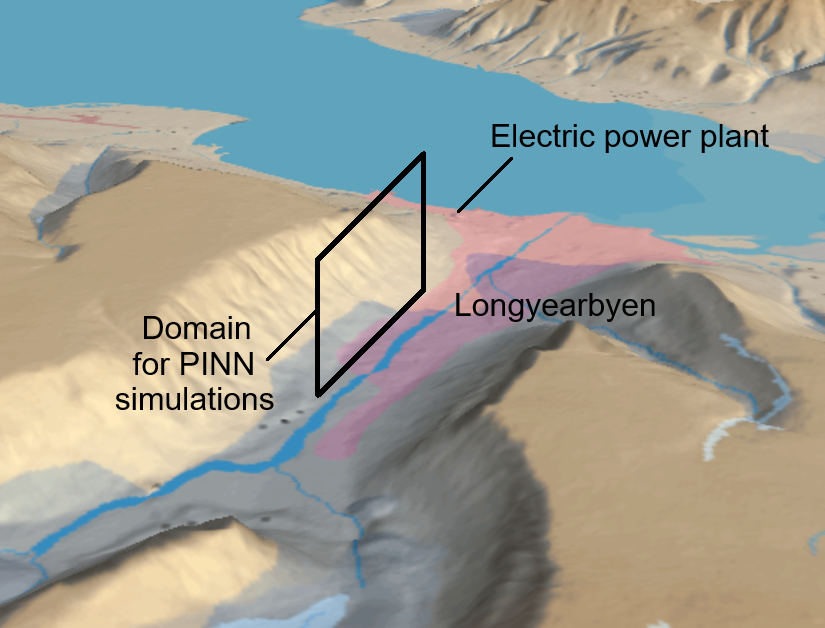}
    \caption{The Town of  Longyearbyen at Spitsbergen, the location of the Diesel power plant, with the computational domain for PINN simulations of thermal inversion (black rhombus.}
    \label{fig:topodomain}
\end{figure}

The subject of our article is the presentation of two computational methods related to simulations of pollution propagation using advection-diffusion equations. The first approach's novelty lies in using the original model of graph grammars, a set of rules describing transformations of the graph describing the computational grid. Based on the generated computational grid, the simulation of pollution propagation is carried out by the finite element method, using a meshless solver and assembling the matrices for the iterative solver from local element matrices. The novelty of the second approach lies in the original implementation of the Physics Informed Neural Networks method in Google Colab, coupled with the results of the finite element method solver, enabling the training of a neural network predicting pollution propagation in the long-term range based on the phenomenon of thermal inversion.
In our simulations, we apply the concepts to a real-world scenario-pollution propagation in the city of Longyearbyen, the capital of Spitzbergen, where air pollution is generated by a coal-burning power plant. This practical application underscores the relevance and importance of our research.

For the simulations of the pollution propagation from the electric power plant, we employ a finite element solver for the advection-diffusion equation \cite{c1} stabilized with the Streamlined-Upwind-Petrov-Galerkin (SUPG) method \cite{c2}.

The computational mesh for the simulations comes from a new graph grammar-based mesh generator, implemented in Julia, for a sequence of mesh refinements built with tetrahedral finite elements. 
Our graph grammar model expresses the three-dimensional version of the longest-edge refinement algorithm.
The longest-edge refinement algorithm has been initially proposed for two-dimensional grids by Cecilia Rivara~\cite{c3,c4}.
The graph grammar-based mesh refinements for two-dimensional grids have been employed and discussed in~\cite{c1}, and in 
\cite{c5,c6,c7,c8} with the hanging nodes version.


%


The topography of the Longyearbyen area has been built using the Global Multi-Resolution Topography (GMRT) synthesis\footnote{\url{https://www.gmrt.org/}}, i.e., a multi-resolution compilation of edited multibeam sonar data collected by scientists and institutions worldwide \cite{c9}.



Finally, to enhance the modeling of the thermal inversion phenomena, we implemented and applied the Physics Informed Neural Networks (PINN) approach~\cite{c10}. 
The PINN simulations of the thermal inversion presented in this paper concern the two-dimensional domain, defined along the valley where the town of Longyearbyen is located (see Figure~\ref{fig:topodomain}). 

The extraordinary success of Deep Learning (DL) algorithms in various scientific fields \cite{c11,c12,c13} over the last decade has recently led to the exploration of the possible applications of (deep) neural networks (NN) for solving partial differential equations (PDEs). The exponential growth of interest in these techniques started with the PINN (\cite{c14}). PINNs have been successfully applied to solve a wide range of problems, from fluid mechanics \cite{c15,c16}, in particular Navier-Stokes equations \cite{c17,c18,c19}, wave propagation \cite{c20,c21}, phase-field modeling \cite{c22}, biomechanics \cite{c23,c24}, quantum mechanics \cite{c25}, electrical engineering \cite{c26}, problems with point singularities \cite{c27}, uncertainty qualification \cite{c28}, dynamic systems \cite{c29,c30}, or inverse problems~\cite{c31,c32, c33}, among others.
In this paper, we use the PINN approach to model the thermal inversion phenomena in the town of Longyearbyen at Spitsbergen.

The rest of the paper is structured as follows:  Section~\ref{sec:gg}, we define
the three-dimensional graph-grammar model expressing Rivara's longest edge refinement algorithm.

Section~\ref{sec:polpropsim} discusses the computational model developed to simulate pollution propagation from the electric power plant chimney. In particular, in Section~\ref{sec:mesh}, we describe the mesh generation process using the graph-grammar code; Section~\ref{sec:strongform} describes the strong form, whereas Section~\ref{sec:weakform} introduces the corresponding weak form of the advection-diffusion equation. Section~\ref{sec:supg} presents the stabilized SUPG formulation, and finally, Section~\ref{sec:numres} summarizes numerical experiments and obtained results.

In  Section~\ref{sec:pinn}, we introduce the PINN concept for simulations of the thermal inversion phenomenon, where we cover, in particular, the loss function employed and the sketch of the applied training algorithm. This is followed by a presentation of the software implementation of the PINN in Section~\ref{sec:code} and its numerical results in Sections~\ref{sec:summer} and ~\ref{sec:winter}. Finally, we conclude the paper in Section~\ref{sec:conclus}.


\section{Graph grammar for the longest edge refinements of three-dimensional tetrahedral elements}
\label{sec:gg}

Mesh refinement lies in subdividing an element of a mesh to obtain a finer mesh. During mesh refinement, the original nodes are not removed, and the topology of the original mesh is preserved.
This is different from re-meshing the domain with smaller-sized elements.
\begin{figure}[!ht]
\centering
\includegraphics[width=.9\textwidth]{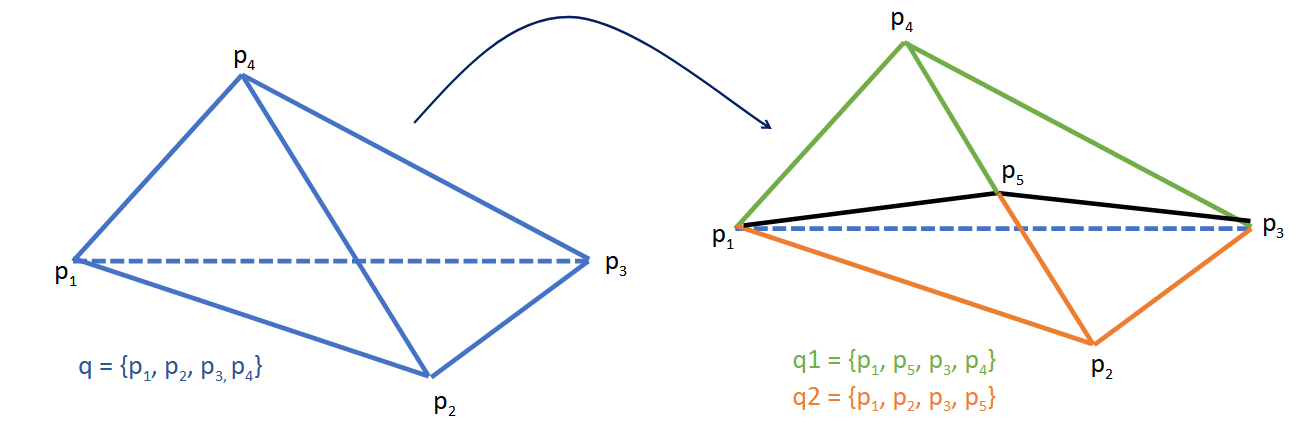}
 \caption{Longest edge refinement for 3D tetrahedral elements.}
\label{fig:LE3D}
\end{figure}
\begin{figure}[!ht]
\centering
 \includegraphics[width=.8\textwidth]{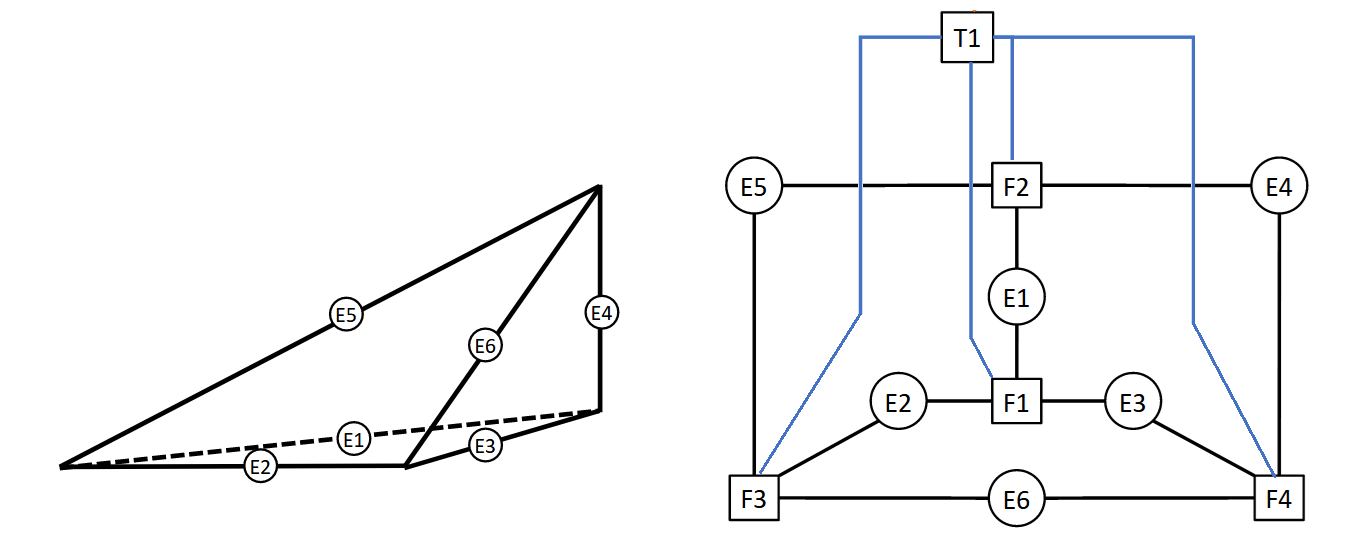}
 \caption{Graph representation of a 3D tetrahedron.}
\label{fig:edges}
\end{figure}
\begin{figure}[!ht]
\centering
\includegraphics[width=.7\textwidth]{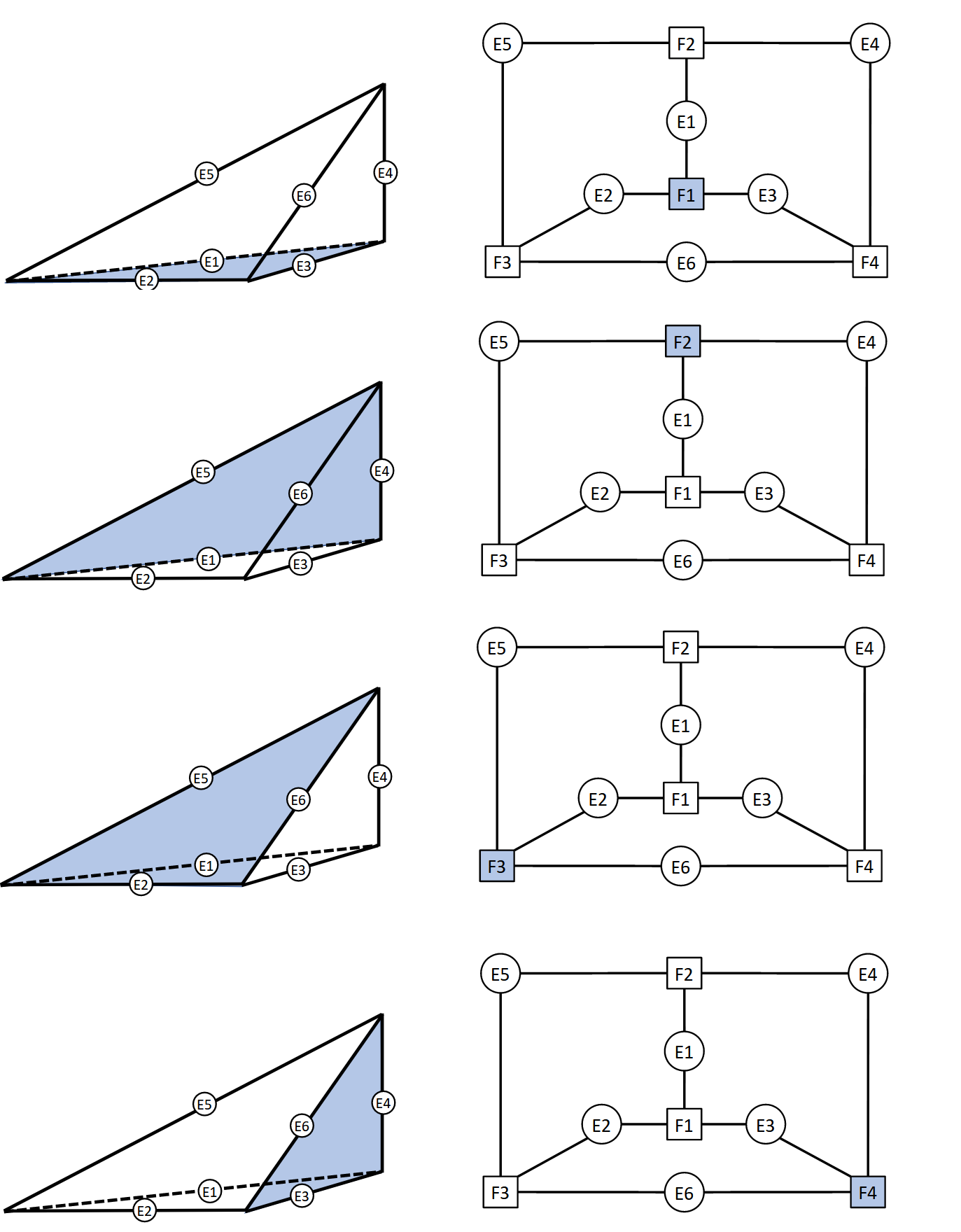}
 \caption{Mapping the tetrahedral faces into the vertices of a graph.}
\label{fig:faces}
\end{figure}
\begin{figure}[!ht] 
\centering
\includegraphics[width=0.4\textwidth]{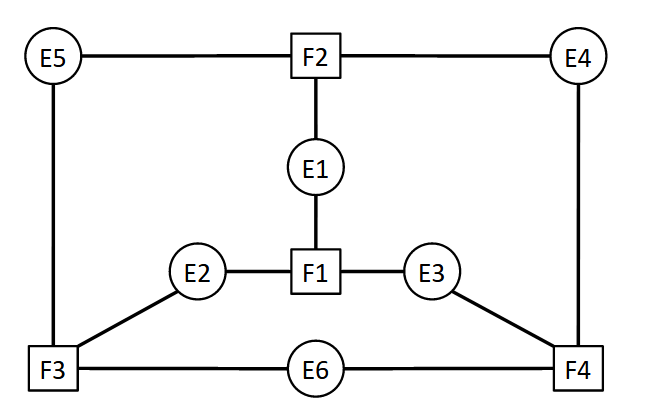}
 \caption{Left-hand side of the graph-grammar production {\bf P1}.}
\label{fig:p1}
\end{figure}
\begin{figure}[!ht]
\centering
\includegraphics[width=.85\textwidth]{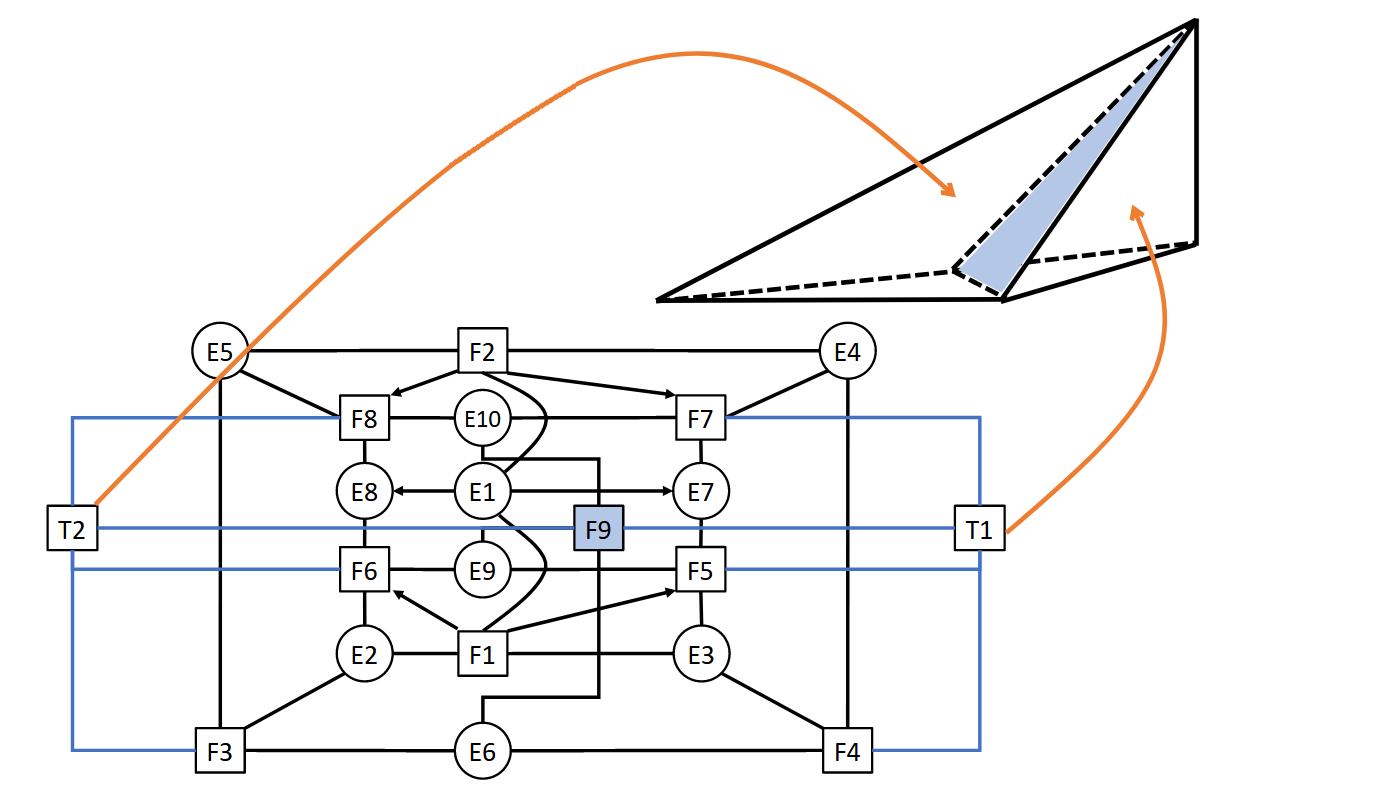} 
 \caption{Right-hand side of the graph-grammar production {\bf P1}, {\bf P2}, {\bf P3}, and {\bf P4}.}
\label{fig:RHS}
\end{figure}
\begin{figure}[!ht] 
\centering
\includegraphics[width=0.48\textwidth]{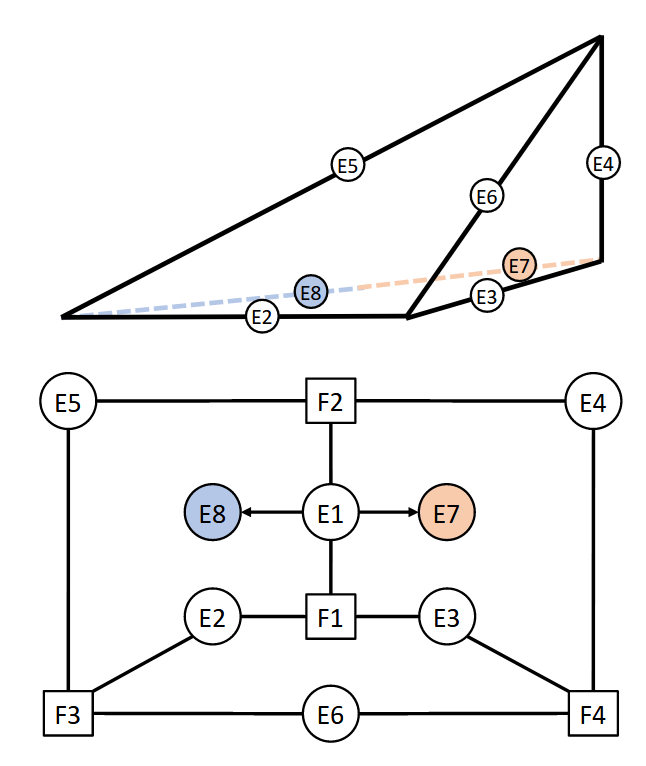}\includegraphics[width=0.48\textwidth]{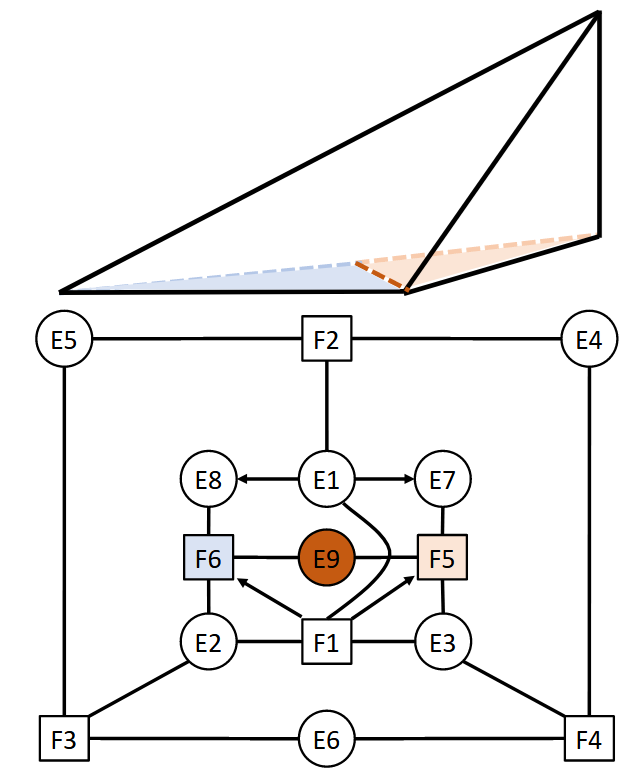} 
 \caption{Graph-grammar production {\bf P2} (left panel) and {\bf P3} (right panel).}
\label{fig:p2}
\end{figure}
\begin{figure}[!ht]
\centering
\includegraphics[width=0.45\textwidth]{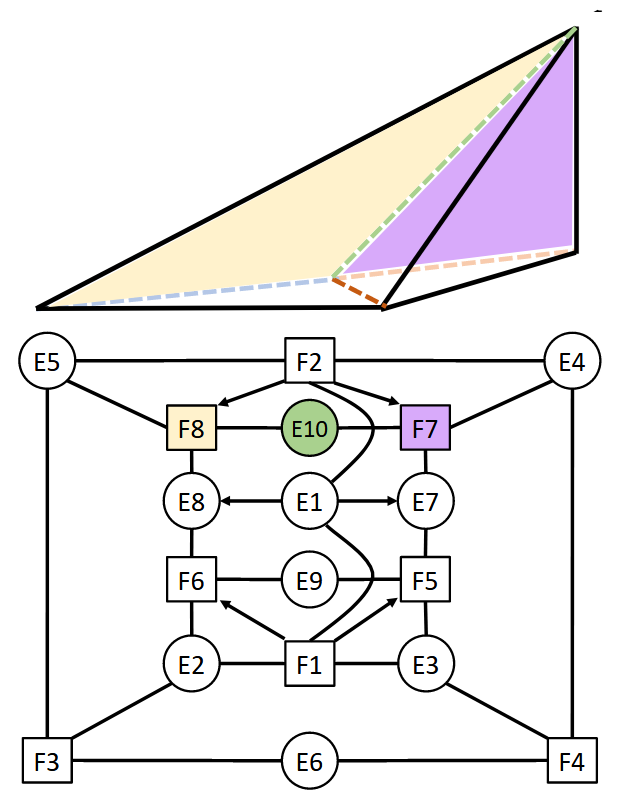} 
 \caption{Graph-grammar production {\bf P4}.}
\label{fig:p4}
\end{figure}
\begin{figure}[!ht]
\centering
    \includegraphics[width=0.32\textwidth]{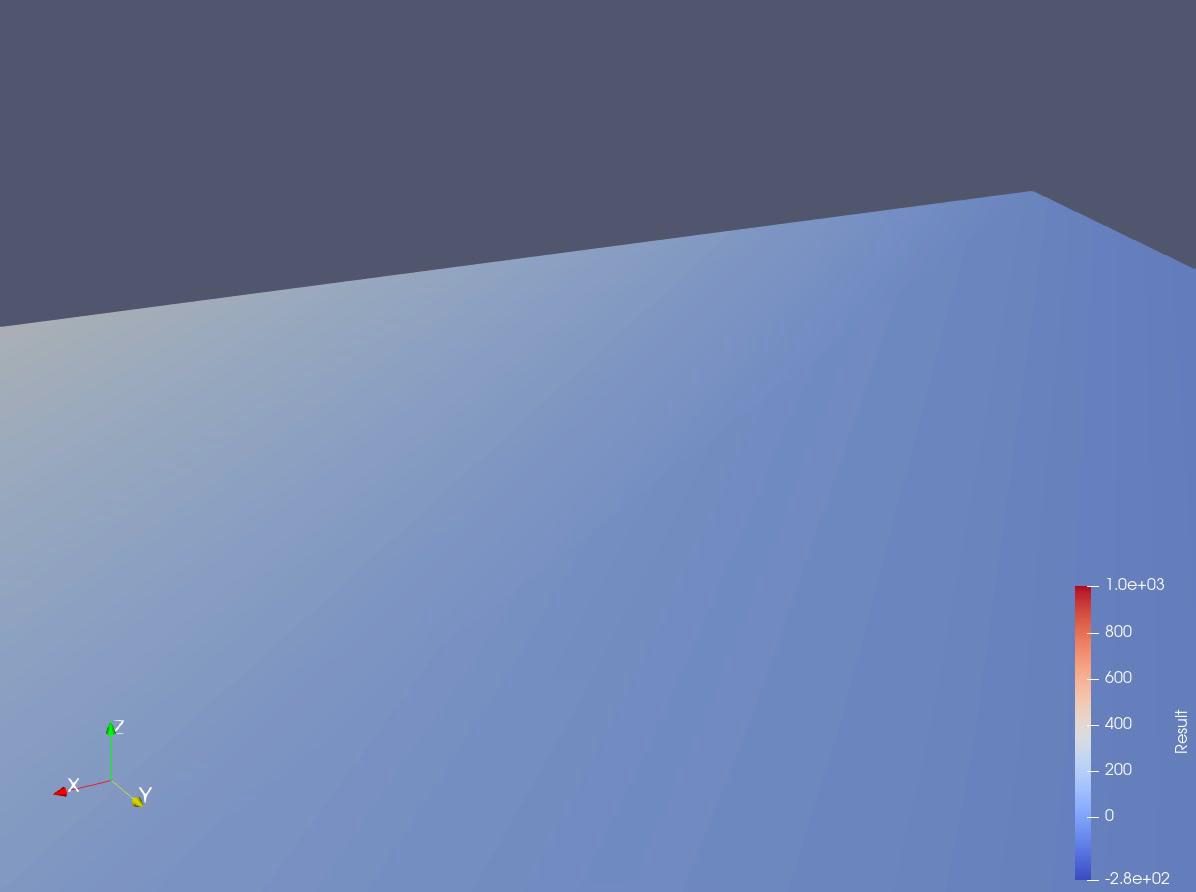}
    \includegraphics[width=0.32\textwidth]{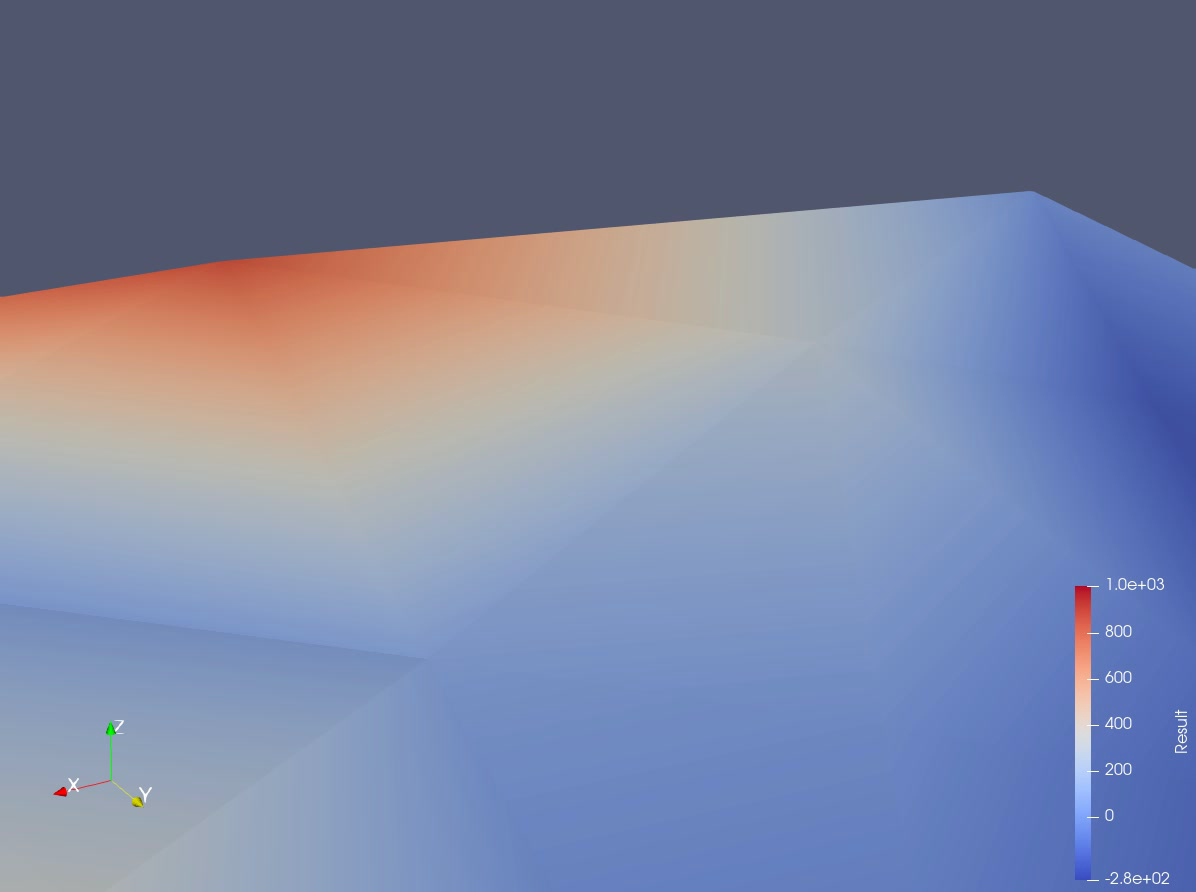}
    \includegraphics[width=0.32\textwidth]{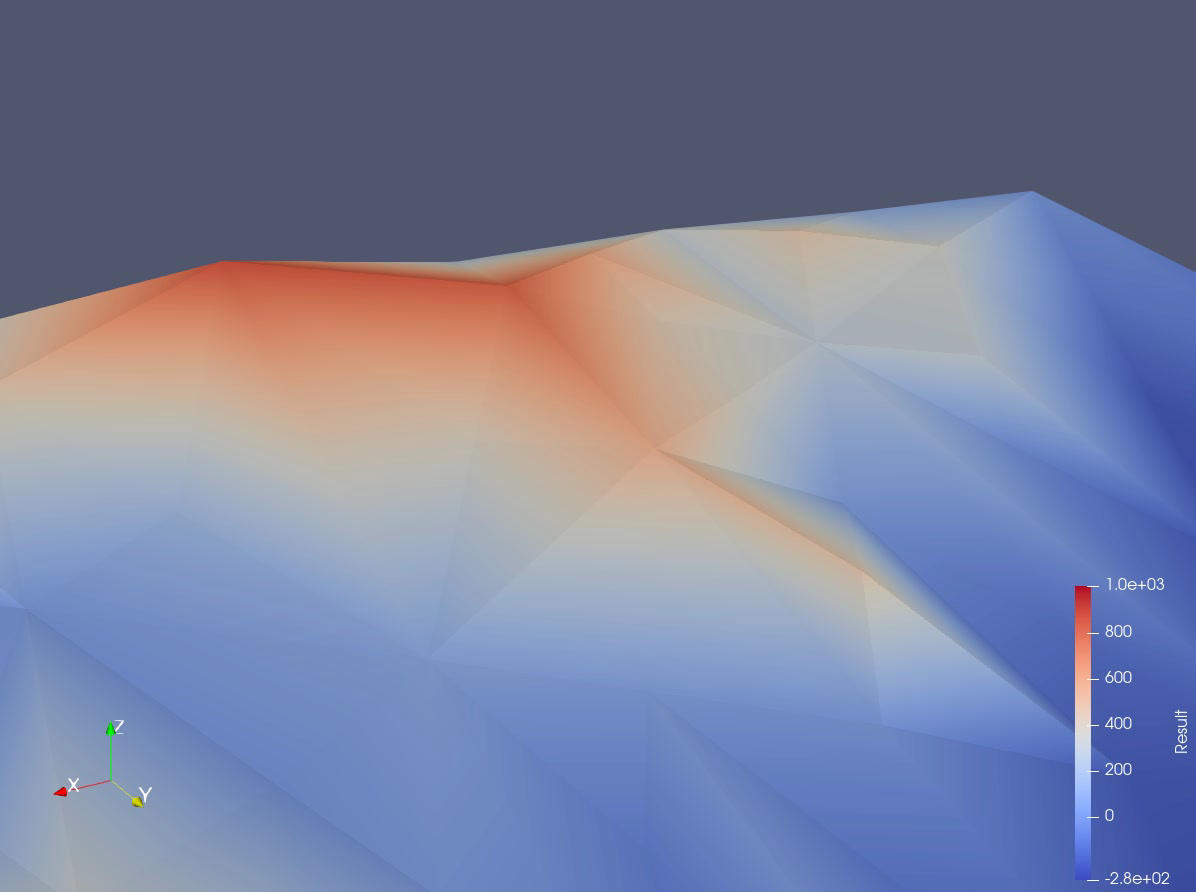}\\
    \includegraphics[width=0.32\textwidth]{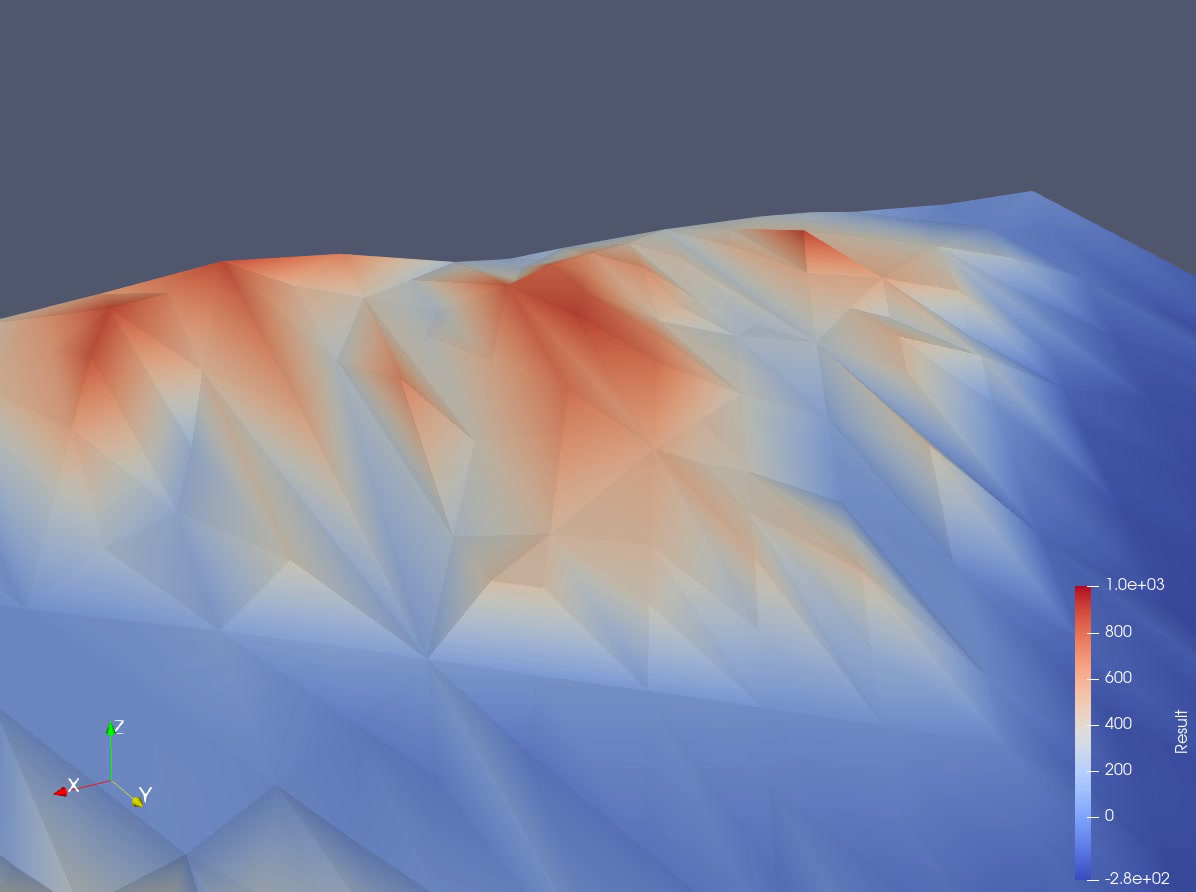}
    \includegraphics[width=0.32\textwidth]{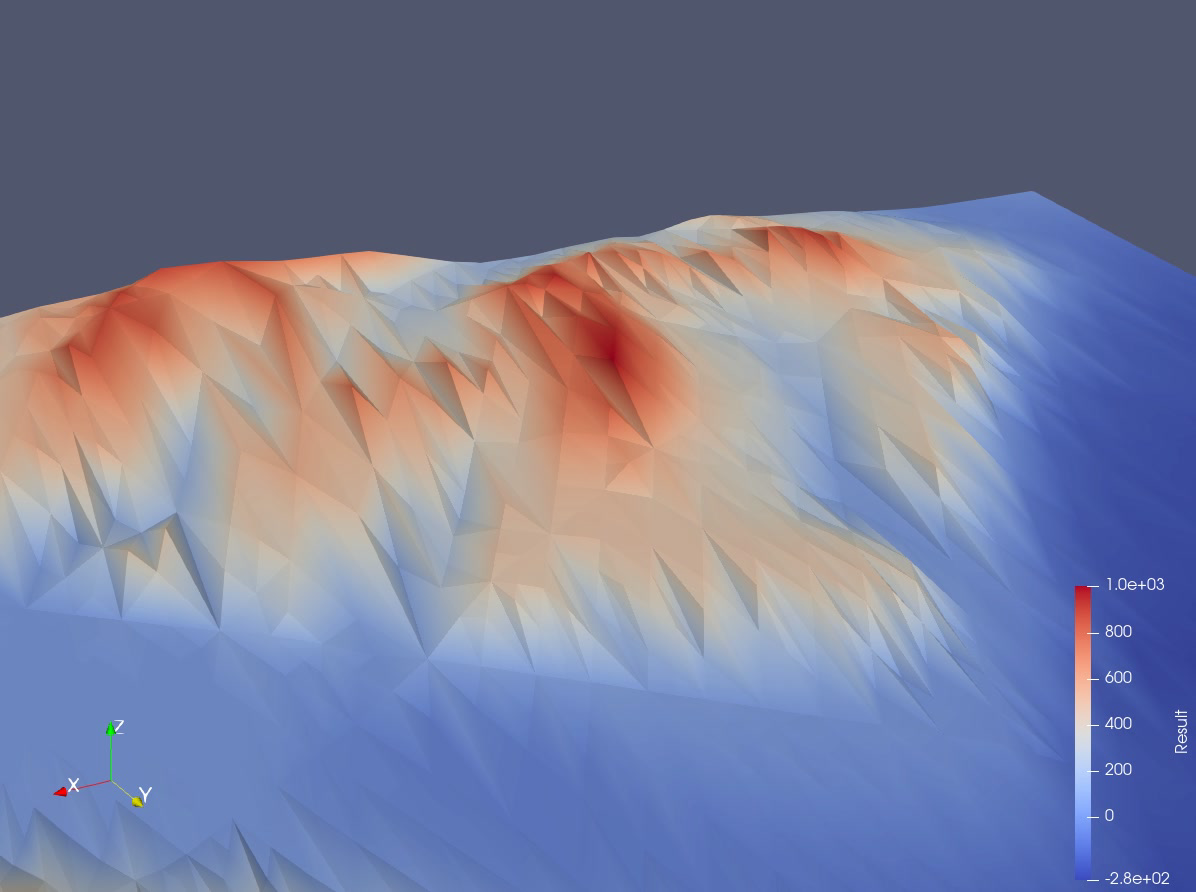}
    \includegraphics[width=0.32\textwidth]{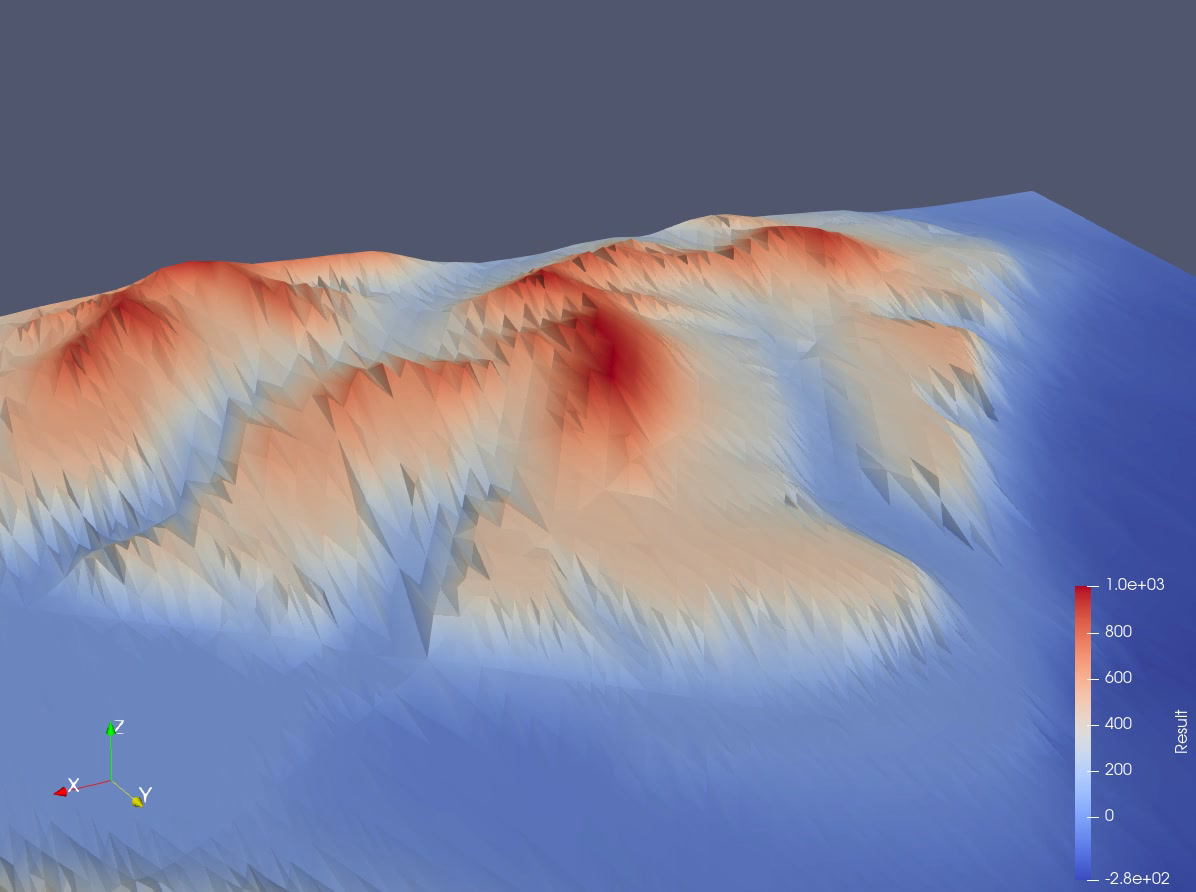}\\
    \includegraphics[width=0.32\textwidth]{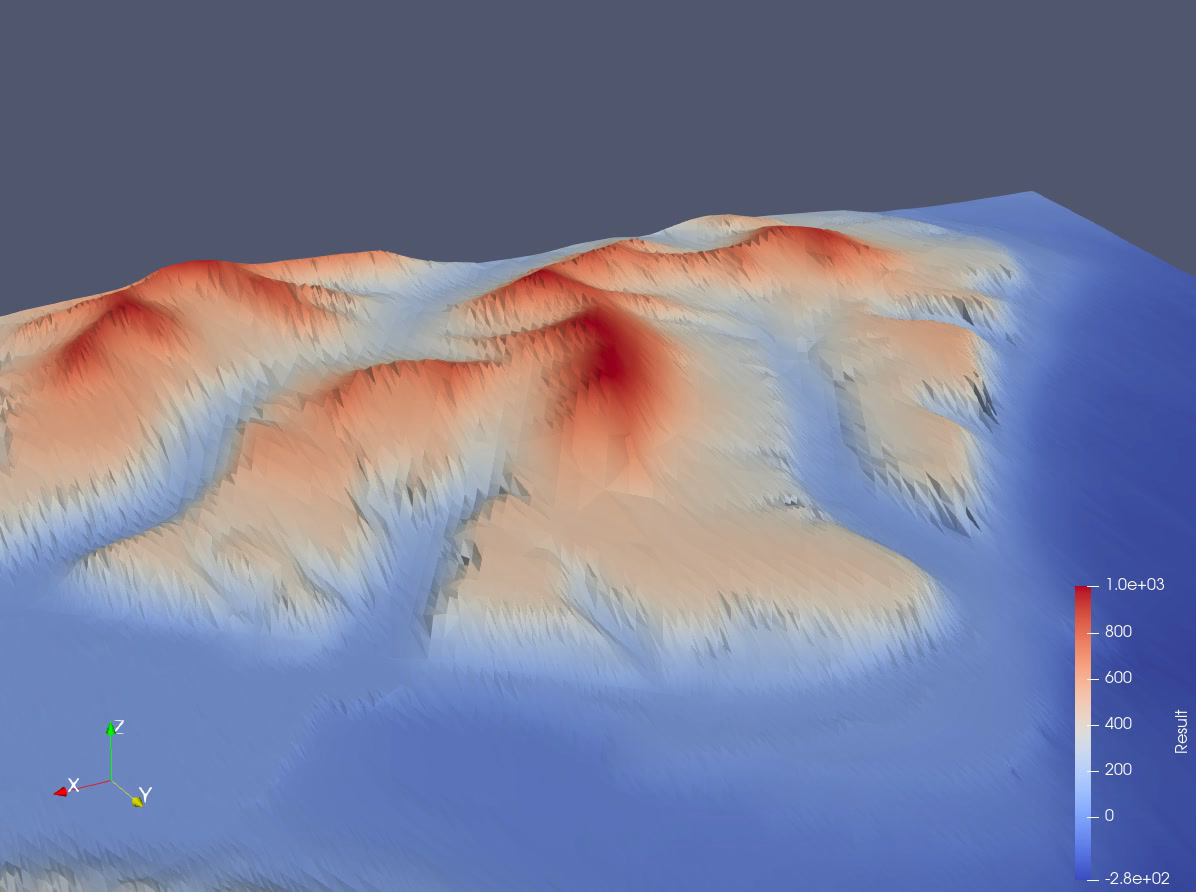}
    \includegraphics[width=0.32\textwidth]{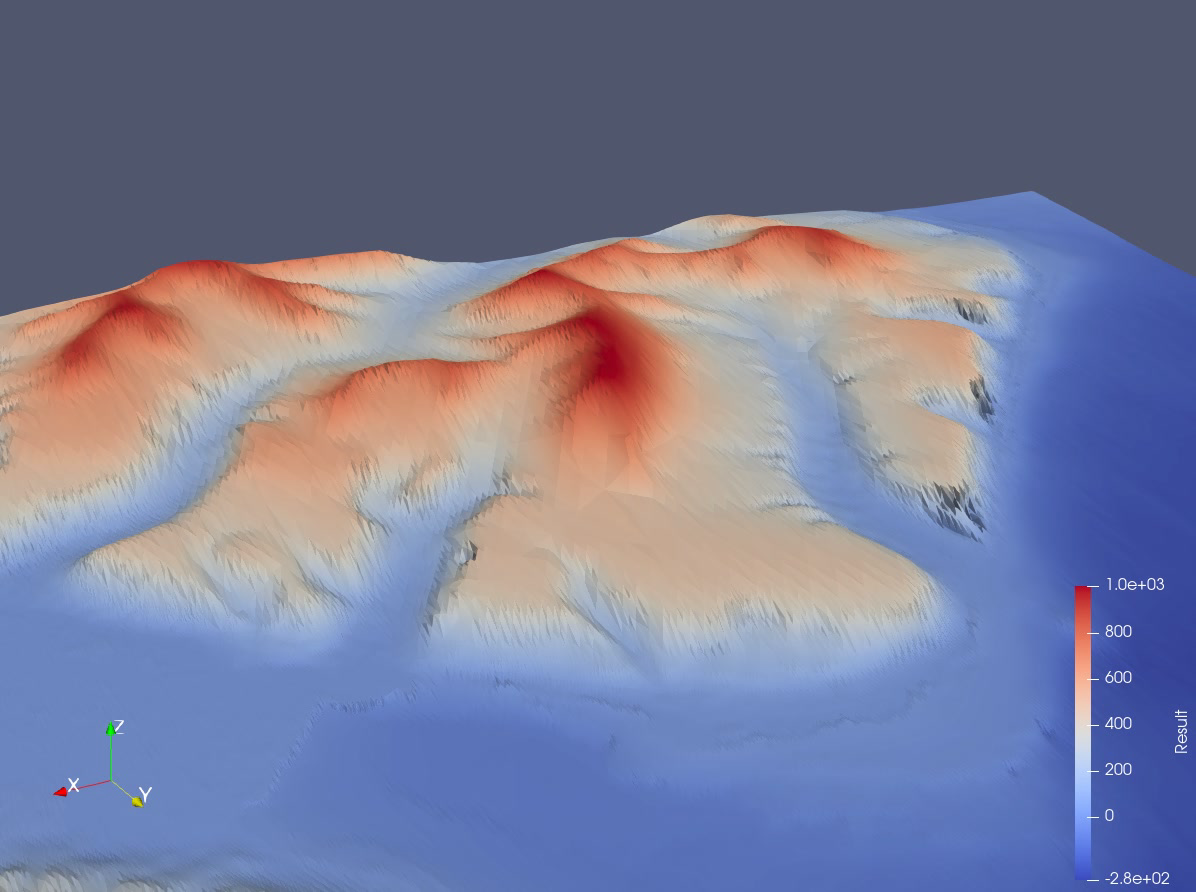}
    \includegraphics[width=0.32\textwidth]{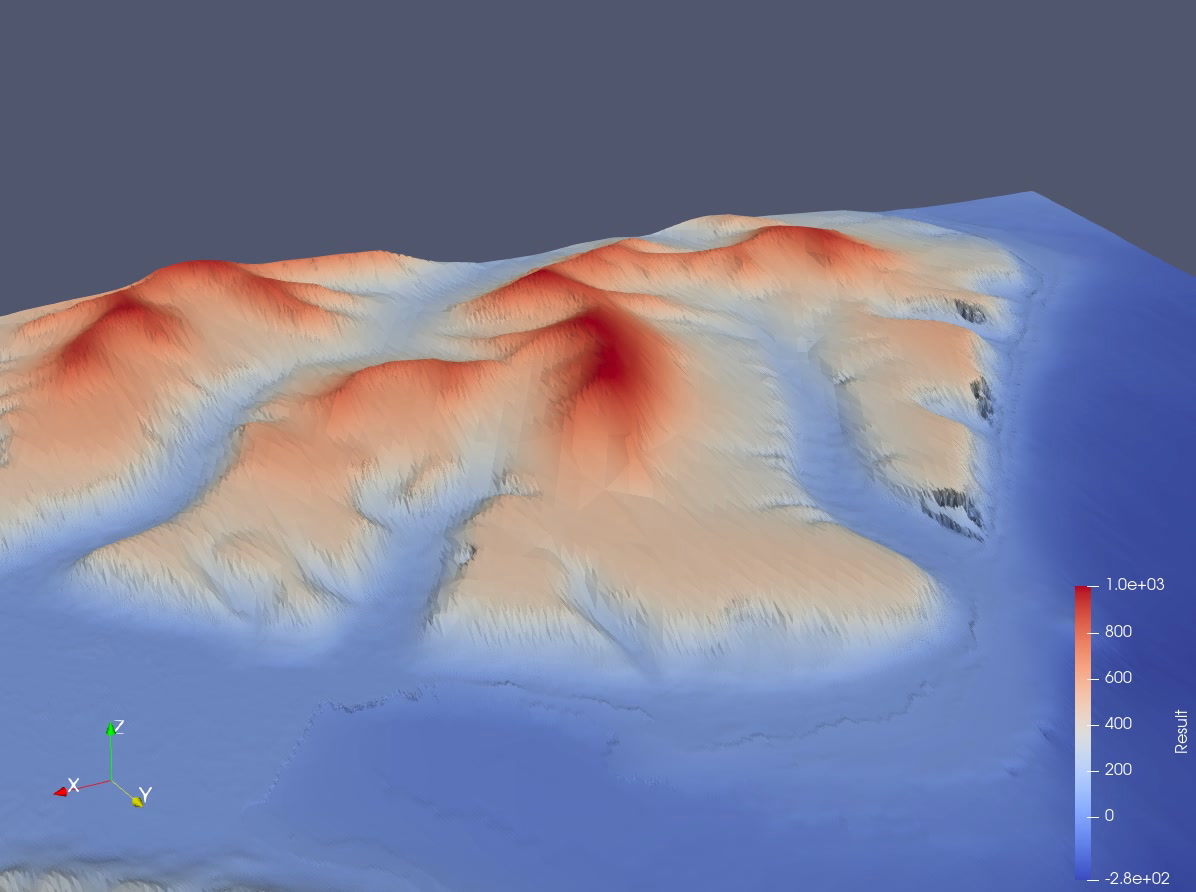}
 \caption{A sequence of mesh refinements performed by the longest-edge refinement algorithm to generate the topography of Spitsbergen.}
\label{fig:sequence}
\end{figure}
Longest-edge refinement (see Figure~\ref{fig:LE3D}) can be expressed mathematically as the bisection of a simplex: 
\begin{equation}
\left(q=\{p_1, p_2,\cdots , p_n, p_{n+1}\}\right) \in {\cal R}^n.
\end{equation}
If the distance between $p_k$ and $p_m$ is the maximum distance of the simplex, then a new point is created such that $p = (p_k+p_m)/2$ and the two new simplices are created such as:
\begin{eqnarray}
q_1 = \{p_1, p_2,\cdots, p_{k-1}, p, p_{k+1}, \cdots, p_m, \cdots, p_{n+1}\} \\
q_2 = \{p_1, p_2, \cdots, p_k, \cdots, p_{m-1}, p, p_{m+1}, \cdots, p_{n+1}\}
\end{eqnarray}

Mathematically, the longest edge can be determined in any dimension.
Geometrically, the longest edge refinement generates a new point in the middle of the longest edge, generating two new elements.
The three-dimensional tetrahedral element is represented as a graph, with vertices representing edges (see Figure~\ref{fig:edges}) and faces (see Figure~\ref{fig:faces}) of the tetrahedron.
We need four graph grammar productions {\bf P1}, {\bf P2}, {\bf P3}, and {\bf P4}  to express for the tetrahedral mesh refinement.
Each of the productions bisects the tetrahedron into two new tetrahedral elements.
The graph representation of the tetrahedral element has the following attributes:
\begin{itemize}
\item Attributes of vertex $T$ representing the whole tetrahedron:
\begin{itemize}
    \item 
$R$: The triangle is marked to be refined
\end{itemize}
\item Attributes of vertex $E$ representing a single edge:
\begin{itemize}
\item $LE$: The edge is one of the longest-edges
\item $BR$: The edge is broken
\item $AE$: The edge is located on the boundary (1 if is a boundary, 2 if is interior)
\item $x, y, z$: Coordinates of the edge (middle point)
\item $IP$: Pointer to the initial point
\item $FP$: Pointer to the final point
\end{itemize}
\item Attributes of vertex $F$ representing a single face:
\begin{itemize}
\item $BRF$: The face is broken
\end{itemize}
\end{itemize}

In the following subsections, we focus on the computational tools developed for the numerical simulations of pollution propagation in Longyearbyen.
In particular,  we introduce a novel graph-grammar model for generating the computational mesh employed for the simulations.

\subsection{Graph-grammar production P1}

The first graph-grammar production denotes the case when the tetrahedral has no broken edges. The production's left-hand side is denoted in Figure~\ref{fig:p1}.
The right-hand side for the graph-grammar production {\bf P1} as well as for all the other productions {\bf P2}, {\bf P3}, and {\bf P4} are presented in Figure~\ref{fig:RHS}.

We have the following predicates of applicability of the graph-grammar production {\bf P1} (\emph{i.e.}, conditions that must be fulfilled if the graph-grammar production can be executed):

\noindent {\tt (NOT BR1 AND LE1) AND ( R1 OR ANY(BRj) ) AND NOT (BRF1 OR BRF2) AND
NOT ANY(BRj AND LEj) 
AND
NOT ANY(NOT BRj AND LEj AND LESS(E1, Ej))},

\noindent The first component  {\tt (NOT BR1 AND LE1)} checks if the first edge is not broken and if it is the longest edge. 
The second component  {\tt ( R1 OR ANY(BRj)}
checks if the tetrahedron has been marked to be refined or already has some broken edges  (is non-conforming) and, therefore, must be broken.
The third component {\tt NOT (BRF1 OR BRF2)}
checks if face1 or face2 are broken, and edge1 is not broken. In this case, the tetrahedron cannot be bisected by edge1.
The fourth component  {\tt NOT (BRj OR LEj)}
checks if any other edge is broken, and it is also the longest edge. In this case, the longest edge is prioritized to be broken
(edge1 will not be broken; we would rather break edge j).
Finally, the fifth component  {\tt NOT ANY(NOT BRj AND LEj AND LESS(E1, Ej))}
checks if any other non-broken edge is also the longest edge.

\subsection{Graph-grammar production P2}
\label{sec:p2}

The second graph-grammar production denotes the case when there is one broken edge of the tetrahedral, but there are no broken faces. The left-hand side of the production is shown in Figure~\ref{fig:p2}, whereas the right-hand side for the graph-grammar production {\bf P2} is presented in Figure~\ref{fig:RHS}.

The following predicates of applicability of the graph-grammar production {\bf P2} exist:

\noindent {\tt (LE1)
AND
NOT (BRF1 OR BRF2)
AND
NOT ANY(BRj AND LEj AND LESS(E1, Ej))}

The first component  {\tt (LE1))}
checks if edge1 is the longest edge, so it should be broken.
The second component {\tt NOT (BRF1 OR BRF2)}
checks if face1 or face2 are broken. In this case, the production {\bf P3} is the right production to apply.
The third component 
{\tt NOT ANY(BRj AND LEj AND LESS(E1, Ej))}
checks if any other broken edge is also denoted as the longest edge. In this case, we will break edge1 only if it is the longest one.

\subsection{Graph-grammar production P3}
\label{sec:p3}

The third graph-grammar production denotes the case when there is one broken edge of the tetrahedral and one adjacent broken face. The production's left-hand side is shown in Figure~\ref{fig:p2} and the right-hand side in Figure~\ref{fig:RHS}.
The following predicates of applicability of the graph-grammar production {\bf P3} exist:

\noindent {\tt (LE1)
AND
(BRF1 AND NOT BRF2)
AND
NOT ANY(BRj AND LEj AND LESS(E1, Ej))}

The first component of the predicate of applicability {\tt (LE1))}
checks if edge1 is the longest edge, so it should be broken.
The second component {\tt (BRF1 AND NOT BRF2)}
checks if face1 is not broken and face2 is broken. In this case. the right production to apply is the production {\bf P2}.
The third component  
{\tt NOT ANY(BRj AND LEj AND LESS(E1, Ej))}
checks if any other broken edge is also denoted as the longest edge. In this case, we will break edge1 only if it is the longest edge.

\subsection{Graph-grammar production P4}
\label{sec:p4}

The fourth graph-grammar production denotes the case when there is one broken edge of the tetrahedral and two adjacent broken faces. The left-hand side of the graph-grammar production is shown in Figure~\ref{fig:p4} and the right-hand side in Figure~\ref{fig:RHS}.

There are the following predicates of applicability of the graph-grammar production {\bf P4}:

\noindent{\tt (LE1)
AND
NOT ANY(BRj AND LEj AND LESS(E1, Ej))}

The first component of the predicate of applicability {\tt (LE1))} checks if edge1 is the longest edge, and, therefore, it should be broken.
The second component {\tt NOT ANY(BRj AND LEj AND LESS(E1, Ej))} checks if any other broken edge is also denoted as the longest edge. In this case, we must use the comparison operator to ensure that edge1 is the longest.

\subsection{Control diagram for graph grammar}
\label{sec:ctrldiag}

\begin{figure}[!ht]
\centering
\includegraphics[width=0.7\textwidth]{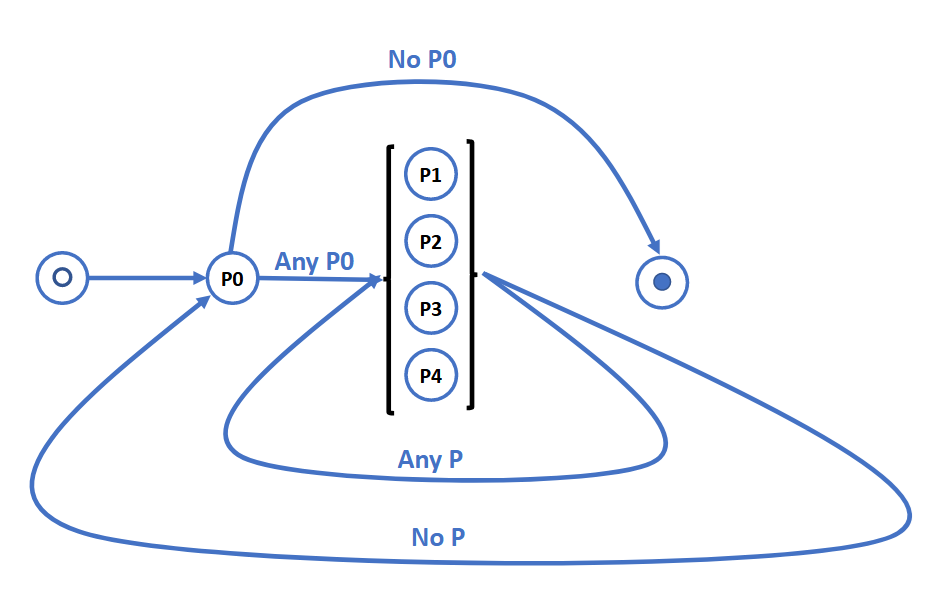} 
 \caption{Diagram controlling the execution of graph-grammar productions {\bf P1}-{\bf P4}.}
\label{fig:diagram}
\end{figure}

The diagram of controlling the execution of the graph grammar production is presented in Figure~\ref{fig:diagram}.

\section{Simulation of the pollution propagation from the electric power plant in Longyearbyen using the Finite Element Method}
\label{sec:polpropsim}

This section describes our finite element method simulations of the advection-diffusion model of pollution propagation from a power plant chimney. 

\subsection{Mesh generation for Spitsbergen topography}
\label{sec:mesh}

We employ the graph grammar described in Section \ref{sec:gg} for the generation of the computational mesh with triangular elements covering the topography of the Longyearbyen area, based on the GMRT data \cite{c9}.
An exemplary sequence of generated meshes is presented in Figure~\ref{fig:sequence}.
In this Figure, we plot the cross-section of the tetrahedral mesh with the approximation of the terrain's topography.
\begin{table}
\centering
\begin{tabular}{cc||cc}
\hline
\textbf{Iteration \#1} & \textbf{Number of nodes}& \textbf{Iteration \#2} & \textbf{Number of nodes}  \\
\hline
1	&	4	&	11	&	8723	\\
2	&	9	&	12	&	16352	\\
3	&	17	&	13	&	29135	\\
4	&	37	&	14	&	49619	\\
5	&	89	&	15	&	83745	\\
6	&	200	&	16	&	144882	\\
7	&	445	&	17	&	258620	\\
8	&	984	&	18	&	440682	\\
9	&	2093	&	19	&	749160	\\
10	&	4355	&	20	&	1572864	\\
\hline
\end{tabular}
\caption{The number of mesh nodes on a generated sequence of triangular element meshes approximating the topography of the Svalbard area, using the 3D longest-edge refinement graph grammar.}
\label{tab:nodes} 
\end{table}
We list the number of generated nodes in Table \ref{tab:nodes}.
The total generation time was equal to 82 seconds.
\begin{figure}[!ht]
\centering
\includegraphics[width=0.9\textwidth]{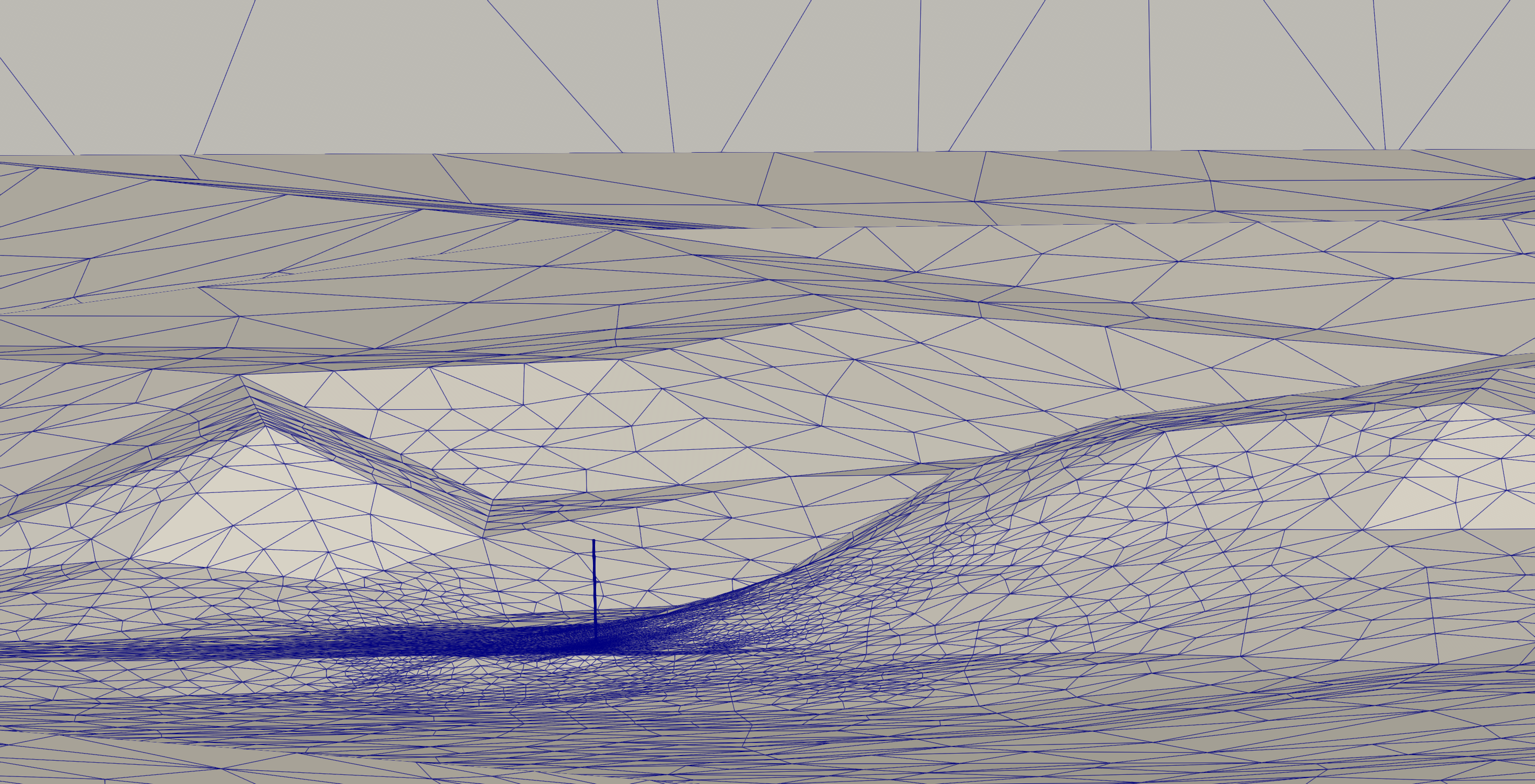}
 \caption{A computational mesh covering the terrain with the chimney.}
\label{fig:mesh1}
\end{figure}
\begin{figure}[!ht]
\centering
\includegraphics[width=0.9\textwidth]{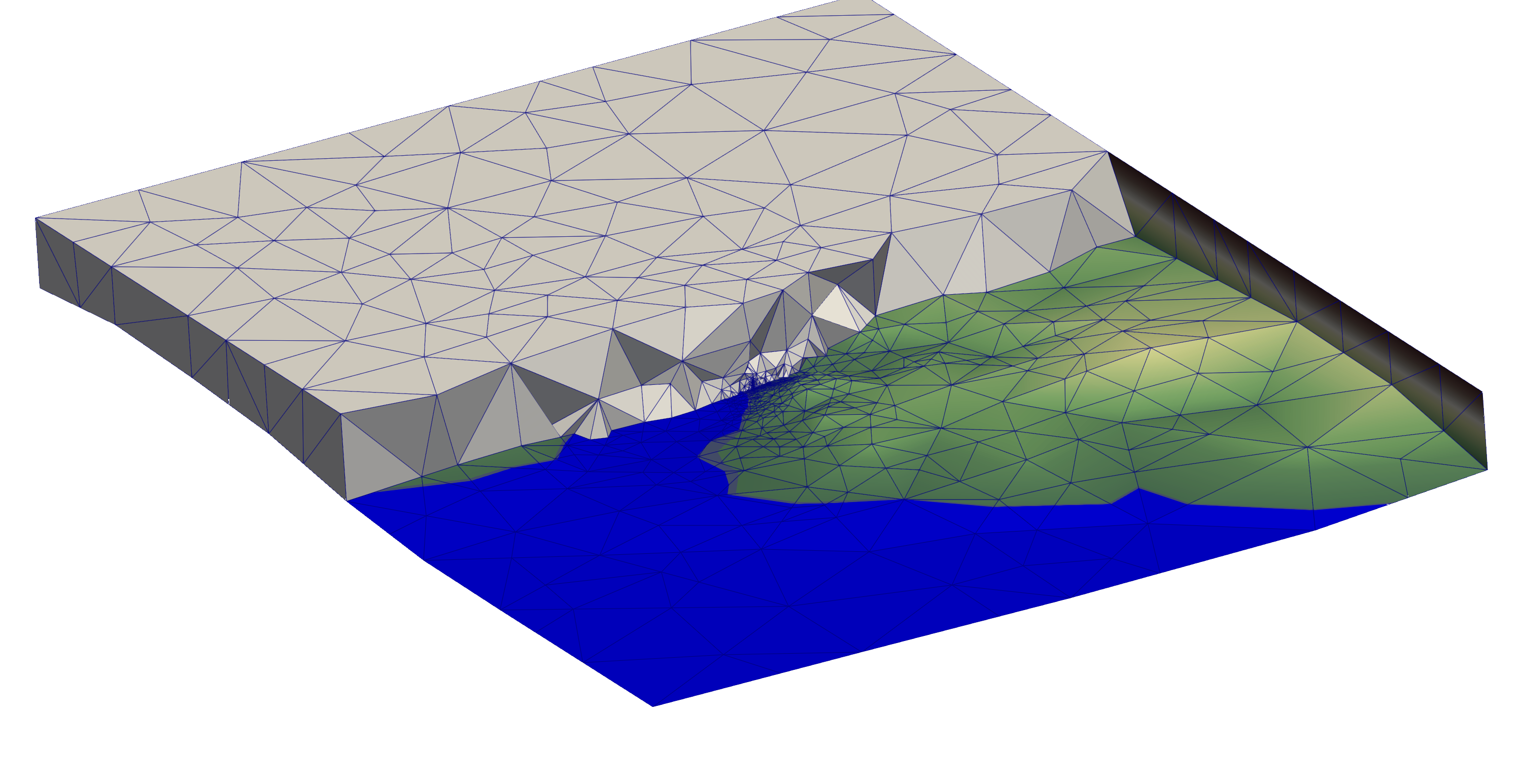} 
 \caption{A cross-section of the 3D computational mesh at the chimney's location.}
\label{fig:mesh2}
\end{figure}
\begin{figure}[!ht]
\centering
\includegraphics[width=0.9\textwidth]{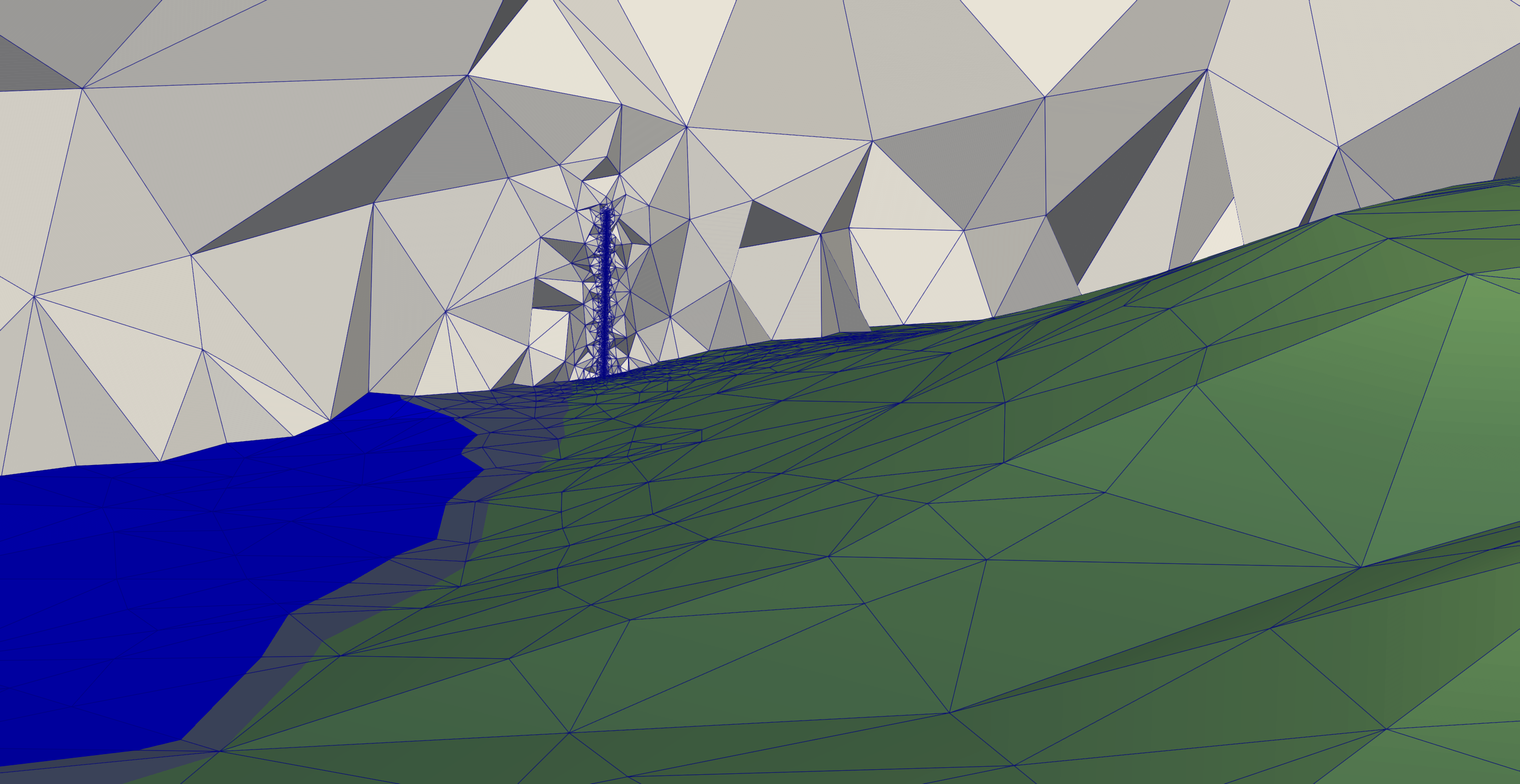}
 \caption{A cross-section of the 3D computational mesh at the chimney's location.}
\label{fig:mesh3}
\end{figure}
The automatically refined mesh from the graph-grammar algorithm has been manually modified to add a chimney representing an electric power plant. For an overview of this mesh, see  Figures~\ref{fig:mesh1},~\ref{fig:mesh2} and~\ref{fig:mesh3}.

\subsection{Strong form of the advection-diffusion-reaction equations}
\label{sec:strongform}

We use the advection-diffusion equation to model the transport of pollutants:
\begin{equation}
\label{eq:confusion}
\frac{\partial u}{\partial t} + \beta \cdot \nabla u - \nabla \cdot \left(\epsilon \nabla u \right)=f,
\end{equation}
where $u(x,y,z,t)$ is the pollutant concentration field; $\beta(x,y,z,t)=(\beta_x(x,y,z,t),\beta_y(x,y,z,t),$ $\beta_z(x,y,z,t))$ the wind velocity vector field, and $\epsilon$  the diffusion coefficient.
%
%
We discretize \eqref{eq:confusion} in time by introducing time steps ${0=t_0<t_1 <t_2<\cdot<t_N=T}$ and the Crank-Nicholson finite difference scheme in time:
\begin{eqnarray}
\frac{u^{t+1}-u^t}{\Delta t} + \beta \cdot \nabla \frac{u^{t+1}+u^t}{2} - 
{\nabla \cdot \left(\epsilon \nabla \frac{u^{t+1}+u^t}{2} \right)+c \frac{u^{t+1}+u^t}{2}}=f^t
\end{eqnarray}

\subsection{Weak form of the advection-diffusion-reaction equations}
\label{sec:weakform}

To apply the finite element method, we introduce the weak formulation of \eqref{eq:confusion} to find $u \in V=H^1(\Omega)$ such that:
\begin{equation}
\frac{u^{t+1}-u^t}{\Delta t} + \frac{b(u^t,v)+b(u^{t+1},v)}{2}=l(v) \quad \forall v\in V
\label{eq:ModelProblem_weak}
\end{equation}
where:
\begin{eqnarray}
b(u,v)=(\beta \cdot \nabla u,v)_{\Omega} -\left( \epsilon \nabla u,\nabla v \right)_{\Omega}+
 (\epsilon n \cdot \nabla u,v)_{\Gamma} + (cu,v)_{\Omega}
\end{eqnarray}
\begin{equation}
l(v)=(f,v)_{\Omega}
\label{eq:ModelProblem_l}
\end{equation}
where we utilize inner product notation: $(u,v)_{\Omega}=\int_{\Omega} uv \text{d}x\text{d}y\text{d}z$, and $(u,v)_{\Gamma}=\int_{\Gamma} uv \text{d}s$ denotes the $L^2$ scalar product on $\Omega$, $\Gamma=\partial \Omega$, and $n=(n_x,n_y,n_z)$ is the vector normal to $\Gamma$.

\subsection{Streamline-Upwind Petrov-Galerkin method}
\label{sec:supg}

For advection-diffusion equations, the standard Bubnov-Galerkin finite element method is known to be numerically unstable for coarse meshes. To make it numerically stable, we apply the Streamline-Upwind Petrov-Galerkin (SUPG) method \cite{c2}. Starting with the Bubnov-Galerkin discretization, we seek for $u_h \in V_h \subset V$ such that:
\begin{eqnarray}
\left(\frac{u_h^{t+1}-u_h^t}{\Delta t},v_h\right) + \frac{b(u_h^t,v_h)+b(u_h^{t+1},v_h)}{2}=l(v_h)  \forall v_h\in V_h\subset V,
\end{eqnarray}
where $V_h$ is span by polynomial functions introduced by the tetrahedral finite elements. The SUPG method modifies then the weak form to stabilize the formulation:
\begin{equation}
\begin{split}
b(u_h^{t+1},v_h) + \sum_K (R(u_h^{t+1}),\tau \beta\cdot \nabla v_h)_K=
l(v_h)+\sum_K (f,\tau \beta\cdot \nabla v_h)_K
 \quad \forall v\in V,
\label{eq:Erikkson_weak_SUPG}
\end{split}
\end{equation}
where $R(u_h^{t+1})=\beta \cdot \nabla u_h^{t+1} +\epsilon \Delta u_h^{t+1}$, and $\tau^{-1}=
\beta \cdot \left(\frac{1}{h^x_K},\frac{1}{h^y_K},\frac{1}{h^z_K} \right) \allowbreak+ 3p^2\epsilon \frac{1}{ {h^x_K}^2+h{^y_K}^2+h{^z_K}^2}$, and $h^x_K,h^y_K$, $h^z_K$ denote three dimensions of an element $K$.
Thus, we have:
\begin{equation}
b_{SUPG}(u_h^{t+1},v_h)=l_{SUPG}(v_h) \quad \forall v_h\in V_h,
\label{eq:SUPG_weak}
\end{equation}
\begin{eqnarray}
\begin{aligned}
b_{SUPG}(u_h^{t+1},v_h)=&{\beta_x}\left(\frac{\partial u_h^{t+1}}{\partial x},v_h\right)_{\Omega}+ {\beta_y}\left(\frac{\partial u_h^{t+1}}{\partial y},v_h\right)_{\Omega}+{\beta_z}\left(\frac{\partial u_h^{t+1}}{\partial z},v_h\right)_{\Omega}+ \nonumber \\
&\epsilon \left(  \frac{\partial u_h^{t+1}}{\partial x}, \frac{\partial v_h}{\partial x}\right)_{\Omega} +\epsilon \left(  \frac{\partial u_h^{t+1} }{\partial y}, \frac{\partial v_h}{\partial y}\right)_{\Omega} + \epsilon \left(  \frac{\partial u_h^{t+1} }{\partial z}, \frac{\partial v_h}{\partial z}\right)_{\Omega}+ \nonumber \\ &(cu_h,v_h)_{\Omega}
-\left(\epsilon\frac{\partial u_h^{t+1}}{\partial x}n_x,v_h\right)_{\Gamma}- \nonumber \\
&\left(\epsilon\frac{\partial u_h^{t+1}}{\partial y}n_y,v_h\right)_{\Gamma}
-\left(\epsilon\frac{\partial u_h^{t+1}}{\partial z}n_z,v_h\right)_{\Gamma}+   \nonumber  \\
&\left(\beta_x\frac{\partial u_h^{t+1}}{\partial x}+\beta_y\frac{\partial u_h^{t+1}}{\partial y}+\beta_z\frac{\partial u_h^{t+1}}{\partial z}+\epsilon \Delta u_h^{t+1}, \right. \nonumber  \\
&\left. \left(\frac{1}{h_x}  + 3\epsilon \frac{p^2}{{h^x_K}^2+{h^y_K}^2}\right)^{-1} 
\beta_x\frac{\partial v_h}{\partial x}+\beta_y\frac{\partial v_h}{\partial y}+\beta_z\frac{\partial v_h}{\partial z}
\right)_{\Omega}
\label{eq:Erikkson_b_SUPG}
\end{aligned}
\end{eqnarray}
\begin{eqnarray}
\begin{aligned}
l_{SUPG}(v_h)&=(f,v_h)_{\Omega}  +\left(f,\left(\frac{1}{h_x}  + 3\epsilon \frac{p^2}{{h^x_K}^2+{h^y_K}^2}\right)^{-1}\right. \left.\left(\beta_x\frac{\partial v_h}{\partial x}+\beta_y\frac{\partial v_h}{\partial y}+\beta_z\frac{\partial v_h}{\partial z}\right)\right)_{\Omega}
. \nonumber
\label{eq:ModelProblem_lsupg}
\end{aligned}
\end{eqnarray}
We incorporate the implicit Crank-Nicholson method into the finite element setup:
\begin{equation}
\begin{split}
&\left( \frac{u^{t+1}-u^t}{\Delta t},w_h\right)_{\Omega}+b_{SUPG}\left(\frac{u_h^t+u_h^{t+1}}{2},v_h\right) = l_{SUPG}(v_h) \quad \forall v_h\in V_h,
\end{split}
\end{equation}
\begin{equation}
\begin{split}
\left(u^{t+1},w_h\right)_{\Omega}+\frac{\Delta t}{2}b_{SUPG}\left(u_h^{t+1},v_h\right)  = \left(u^t,w_h\right)_{\Omega}+\frac{\Delta t}{2}b_{SUPG}\left(u_h^t,v_h\right)+ & l_{SUPG}(v_h) \nonumber \\ &  \forall v_h\in V_h.
\end{split}
\end{equation}
The element matrices and right-hand-side vectors are discretized to obtain
the local systems over each element, with matrices and right-hand-sides:
\begin{eqnarray}
\begin{bmatrix}
({\psi}_1,{\psi}_1) & \cdots & ({\psi}_1,{\psi}_{15}) \\
\vdots & \vdots & \vdots \\
({\psi}_{15},{\psi}_1) & \cdots & ({\psi}_{15},{\psi}_{15})
\end{bmatrix}\begin{bmatrix}
u^t_1 \\
\vdots \\
u^t_{15}
\end{bmatrix} + \nonumber\\ \frac{\Delta t}{2}* 
\begin{bmatrix}
 b^K_{SUPG}({\psi}_1,{\psi}_1) & \cdots & b^K_{SUPG}({\psi}_1,{\psi}_{15}) \\
\vdots & \vdots & \vdots \\
b^K_{SUPG}({\psi}_{15},{\psi}_1) & \cdots & b^K_{SUPG}({\psi}_{15},{\psi}_{15})
\end{bmatrix}
\begin{bmatrix}
u^t_1 \\
\vdots \\
u^t_{15}
\end{bmatrix}=
\begin{bmatrix}
l^K_{SUPG}({\psi}_1) \\
\vdots \\
l^K_{SUPG}({\psi}_{15})
\end{bmatrix}
\label{RHS}
\end{eqnarray}
The resulting local systems are submitted to the matrix-free GMRES iterative solver.

\subsection{Numerical results}
\label{sec:numres}

\begin{figure}[!ht]
\centering
\includegraphics[width=0.8\textwidth]{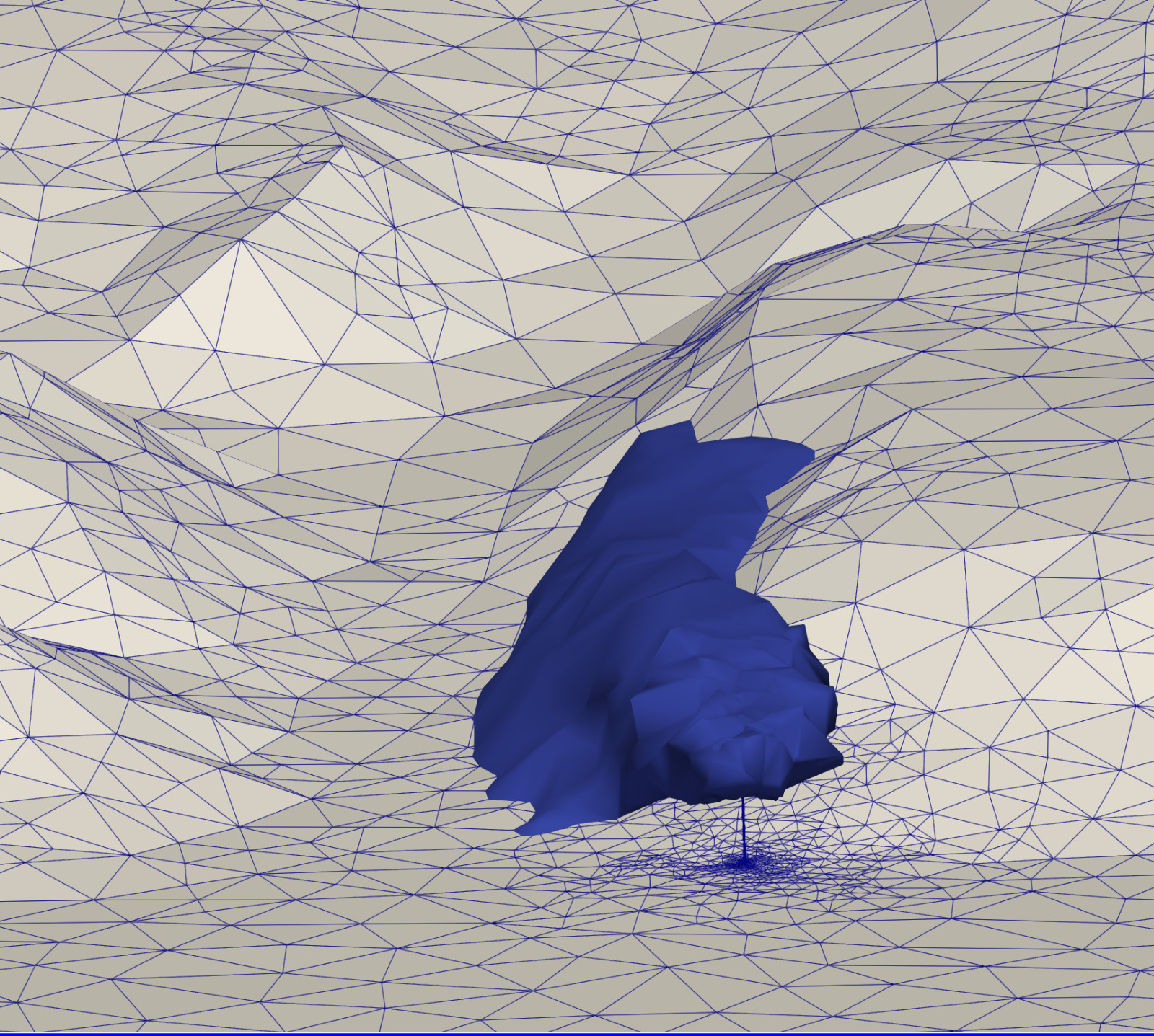}   
\caption{The front view of the smoke propagated from the chimney into the valley after two hours of power plant operation.}
\label{fig:smoke1a}
\end{figure}

\begin{figure}[!ht]
\centering
\includegraphics[width=0.8\textwidth]{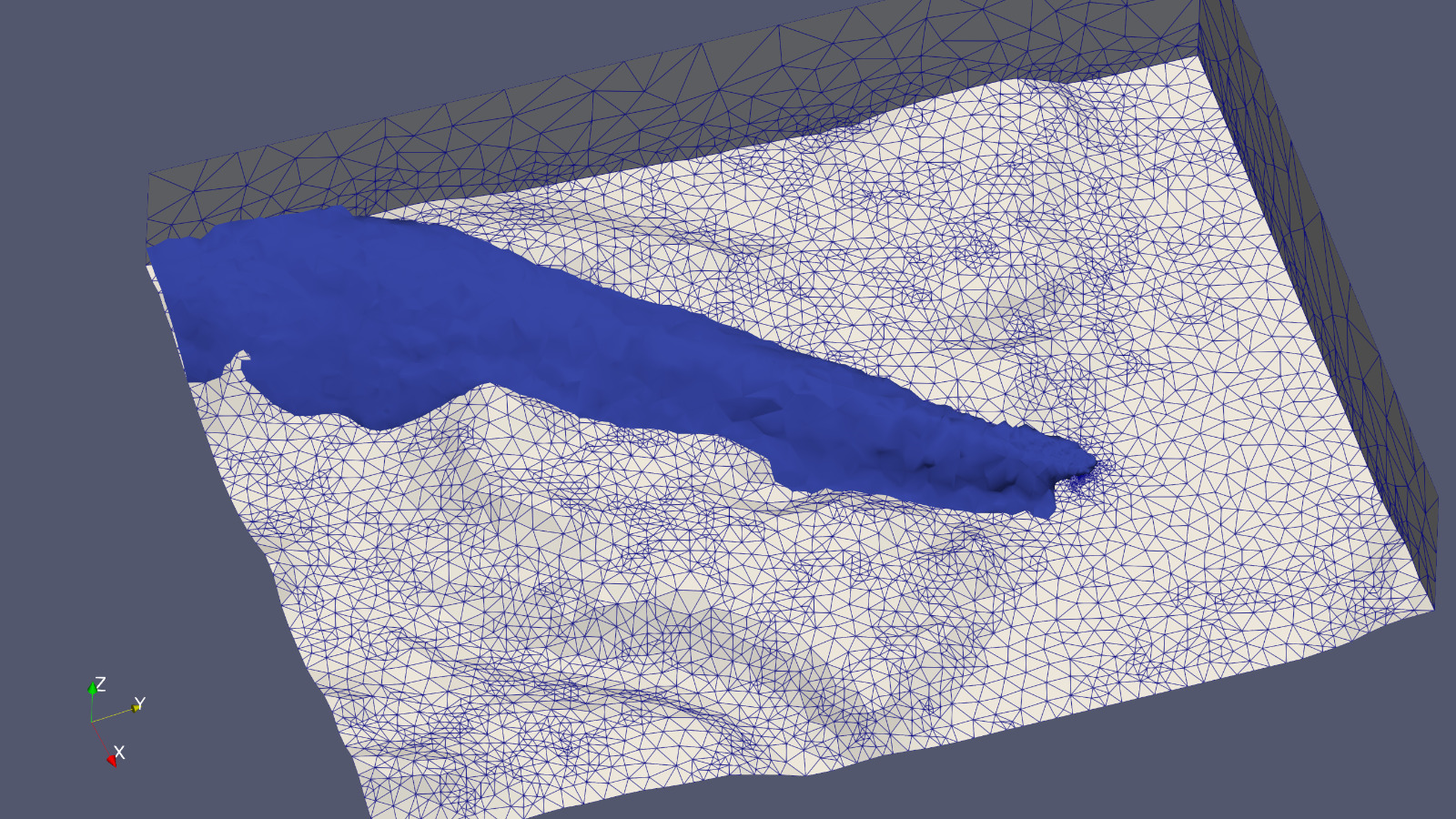}   
\caption{The top view of the smoke propagated from the chimney into the valley after 9 hours of power plant operation.}
\label{fig:smoke1b}
\end{figure}

\begin{figure}[!ht]
\centering
\includegraphics[width=0.8\textwidth]{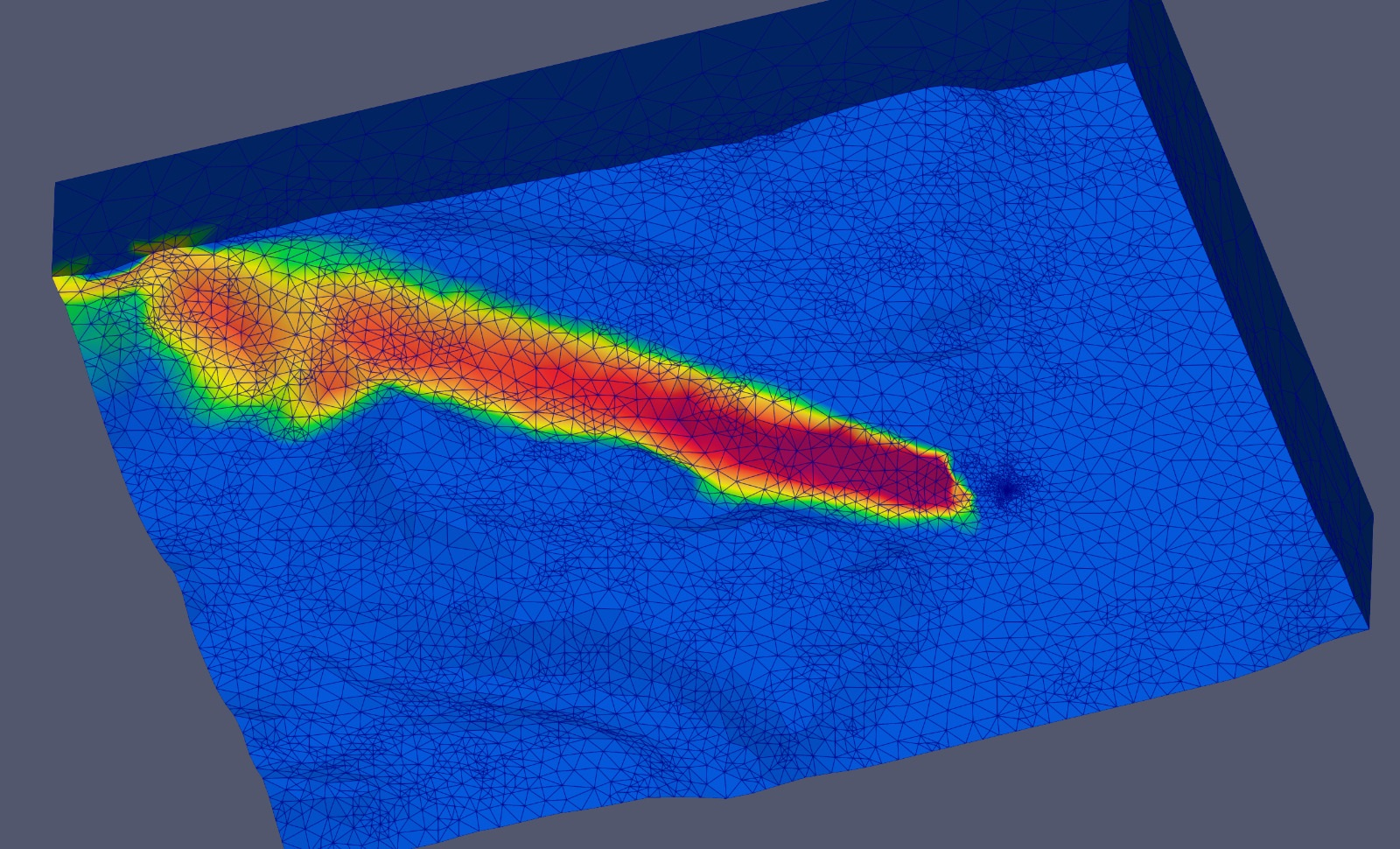}   
\caption{The concentration of the pollution near the ground after 9 hours of working of the electric power plant.}
\label{fig:smoke2}
\end{figure}

We simulated the pollution propagation with the source located on the top of the chimney, assuming the average wind direction and velocity as for the winter season.
As illustrated in Figures \ref{fig:smoke1a}-\ref{fig:smoke2}, the pollution propagates into the valley where Longyearbyen is located.
Nine hours after the chimney starts producing the pollution, the whole valley is filled with pollution.

\section{Simulation of thermal inversion using Physics Informed Neural Networks}
\label{sec:pinn}

\begin{figure}[!ht]
    \centering
    \includegraphics[width=\textwidth]{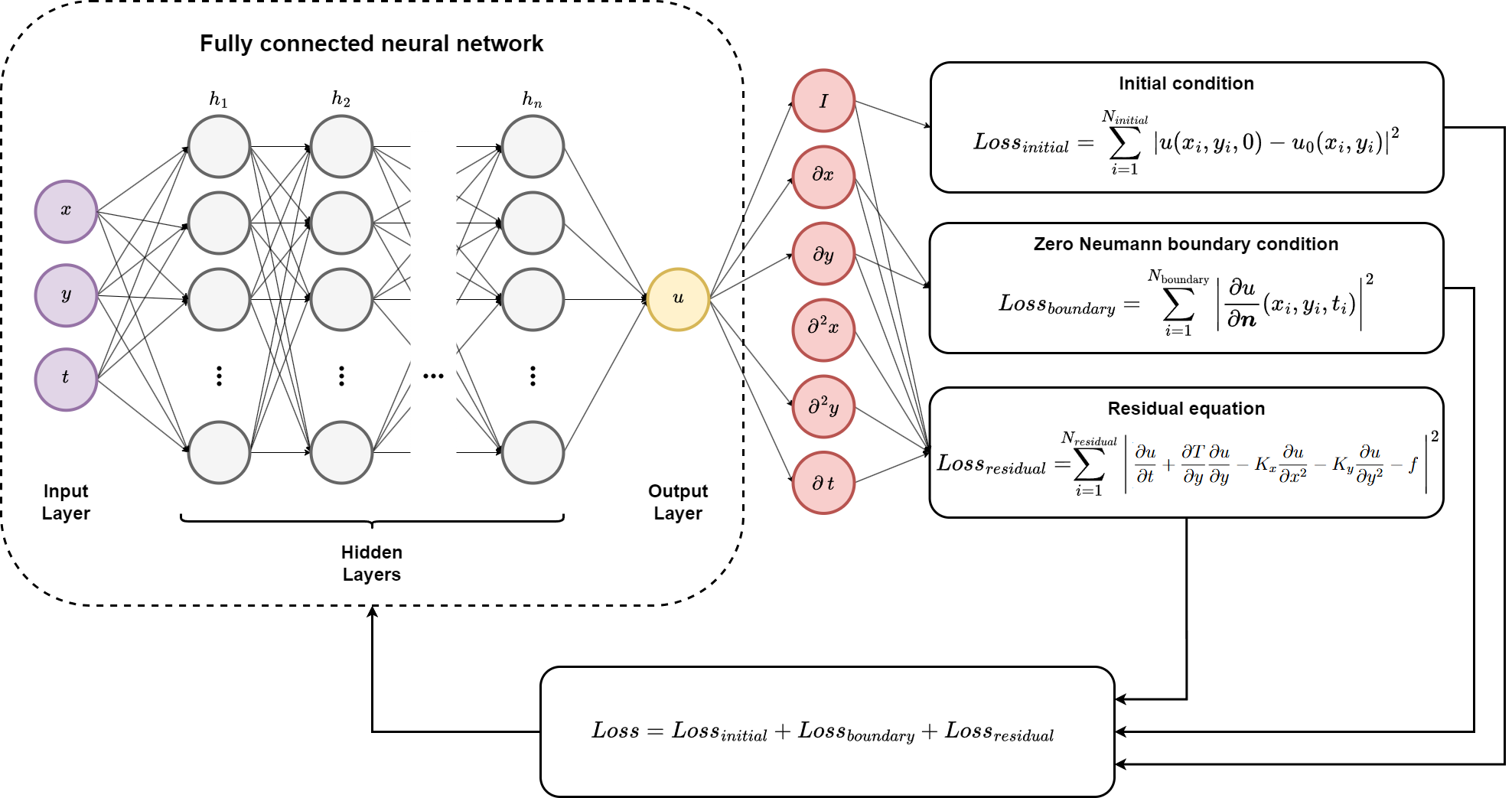}
    \caption{The structure of the Physics Informed Neural Network for modeling of time-dependent advection-diffusion equations.}
    \label{fig:pinn_ad_diagram}
\end{figure}

This section discusses modeling the thermal inversion effect with the Physics Informed Neural Networks \cite{c15,c14}.

\begin{figure}[!ht]
    \centering
    \includegraphics[width=0.48\textwidth]{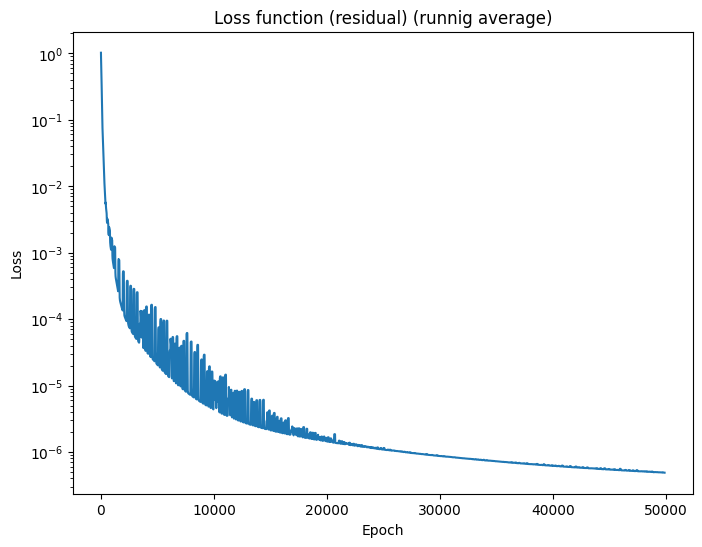}
    \includegraphics[width=0.48\textwidth]{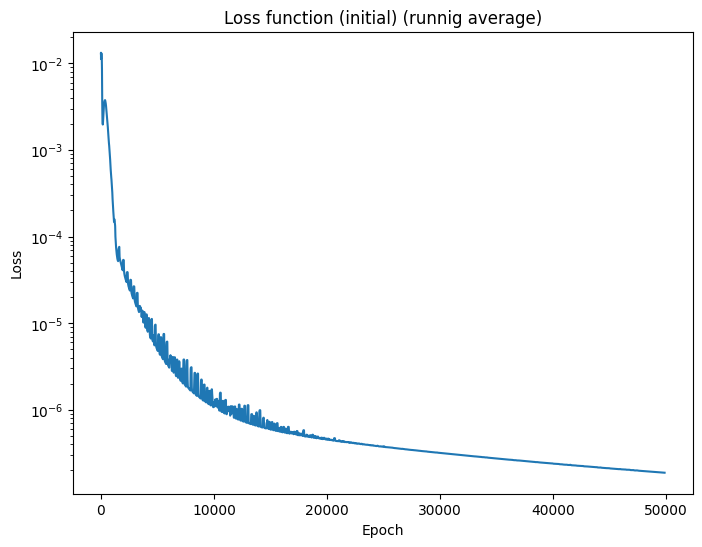}\\
    \includegraphics[width=0.48\textwidth]{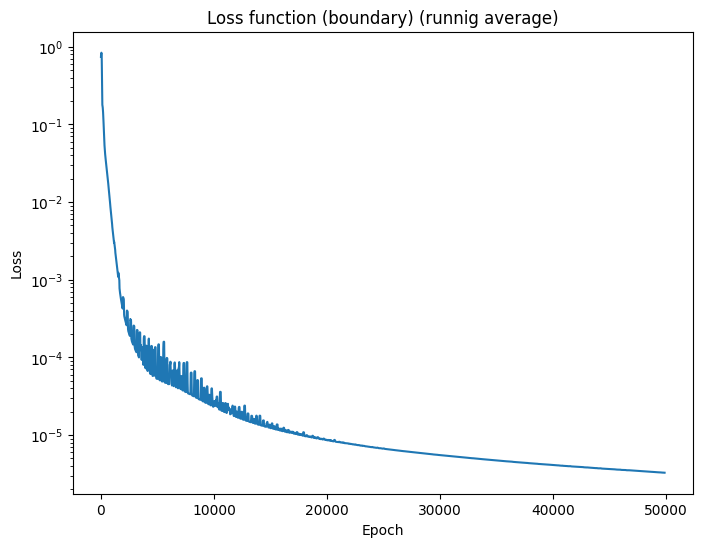}
    \includegraphics[width=0.48\textwidth]{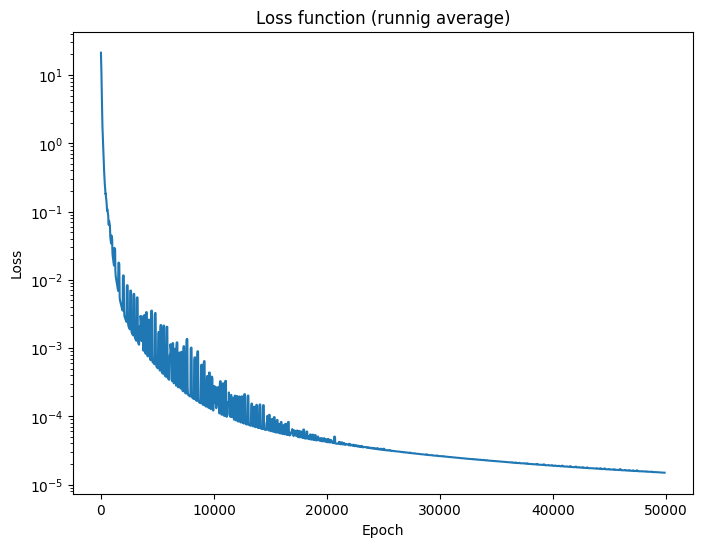}
    \caption{Svalbard summer. The temperature decreases in vertical direction close to the ground. The convergence of residual, initial, boundary, and total loss functions.}
    \label{fig:thermal:losses_summer}
\end{figure}

The PINN code used for thermal inversion simulations are available at

{\tt https://colab.research.google.com/drive/15WDZZV36v2qmzvU\_vq0Ter0RvKxwrYs9 }

and 

{\tt https://colab.research.google.com/drive/1Ta29ihEOX6rWhDozK\_7u89Ev0A3dX8sz }

As the initial state for the simulation, we consider 
the pollution propagated into the valley by the chimney as directed by the wind. This pollution concentration is based on the finite element method solver.
The vertical temperature profile effect is obtained by introducing the advection field as the temperature gradient.

Thermal inversion, also known as temperature inversion, is a meteorological phenomenon where the typical (decreasing with height) temperature gradient of the atmosphere is reversed (increasing with height). Typically, the temperature decreases with altitude, meaning the air is cooler higher up. However, during a thermal inversion, a layer of cooler air becomes trapped near the ground by warmer air above it.
The trapped cold air traps also near the ground pollutants, leading to poor air quality. 
Thermal inversions are more common in valleys, where the topography limits air circulation.
Inversions are also more likely to occur during the winter, especially during clear nights when the ground cools rapidly.

Thermal inversions are a common and significant phenomenon during the Arctic night due to the extreme and prolonged cold conditions that characterize this region. The Arctic night refers to the period during the winter months when the Sun does not rise above the horizon for an extended period, resulting in continuous darkness.
During the thermal inversion phenomenon, the temperature increases with altitude instead of decreasing. Inversions are typical in winter when the low layers of the atmosphere are cooled by a cold surface covered with snow and ice while the higher layers remain warmer.

Thermal inversions during the Arctic night are a natural consequence of the region's extreme and prolonged cold conditions. They result in very stable and cold air near the surface, with warmer air above, and can persist for long periods. 

We model the thermal inversion by introducing the vertical temperature profiles specific to winter and summer seasons in the region of the Town of Longyearbyen at Spitsbergen.

We also assume that the horizontal diffusion coefficient $K_x=0.1$ is stronger than the vertical diffusion coefficient $K_y=0.01$. 
We focus on advection-diffusion equations in the strong form. We seek the pollution concentration field $[0,1]^2\times [0,1] \ni (x,y,t) \rightarrow u(x,y,t) \in {\cal R}$  

\begin{eqnarray}
\frac{\partial u(x,y,t)}{\partial t}+  \left( b(x,y,t) \cdot \nabla \right)u(x,y,t)- \nabla \cdot \left(K \nabla u(x,y,t)\right)  = 0,
  \\ (x,y,t) \in \Omega \times (0,T] \nonumber \\
\frac{\partial u(x,y,t)}{\partial n} = 0, \textrm{  } (x,y,t) \in \partial \Omega \times (0,T] \\
u(x,y,0) = u_0(x,y), \textrm{  } (x,y,t) \in \Omega \times 0 
\end{eqnarray}
This PDE translates into
\begin{eqnarray}
\frac{\partial u(x,y,t)}{\partial t}+  \frac{\partial T(y)}{\partial y}\frac{\partial u(x,y,t)}{\partial y} - 0.1\frac{\partial  u(x,y,t)}{\partial x^2}-0.01\frac{\partial  u(x,y,t)}{\partial y^2}  = 0, \\
  (x,y,t) \in \Omega \times (0,T] \nonumber \\
\frac{\partial u(x,y,t)}{\partial n} = 0, \textrm{  } (x,y,t) \in \partial \Omega \times (0,T] \\
u(x,y,0) = u_0(x,y). \textrm{  } (x,y,t) \in \Omega \times 0 
\end{eqnarray}

Neural networks are composed of interconnected layers of nodes, or neurons, designed to process and learn from data. The architecture of a typical neural network is shown in Figure \ref{fig:pinn_ad_diagram}.
The input layer receives the independent variables of the problem $(x,y,t)$.
The hidden layers are crucial for learning complex patterns. Each neuron in a hidden layer applies a nonlinear $\sigma(x) = \frac{e^x - e^{-x}}{e^x + e^{-x}}$ activation function to a weighted sum of the inputs. Mathematically, the transformation at the \(l\)-th hidden layer is given by:

\[
\mathbf{h}^{(l)} = \sigma \left( \mathbf{W}^{(l)} \mathbf{h}^{(l-1)} + \mathbf{b}^{(l)} \right)
\]

where:
\begin{itemize}
    \item \( \mathbf{h}^{(l)} \) is the output of the \(l\)-th hidden layer,
    \item \( \sigma \) denotes the activation function (e.g., tanh),
    \item \( \mathbf{W}^{(l)} \) represents the weight matrix for the \(l\)-th layer,
    \item \( \mathbf{b}^{(l)} \) is the bias vector for the \(l\)-th layer,
    \item \( \mathbf{h}^{(l-1)} \) is the output of the \((l-1)\)-th layer (or the input layer for \(l=1\)).
\end{itemize}
The output layer generates the final prediction of resulting pollution concentration field $u(x,y,t)$ at point $(x,y,t)$. The output layer is computed as:

\[
u = \mathbf{W}^{(L)} \mathbf{h}^{(L-1)} + \mathbf{b}^{(L)}
\]

where \( L \) denotes the total number of layers in the network.

We define the loss function as the residual of the PDE:

\begin{eqnarray}
L_{residual}(x,y,t) =
\left(
\frac{\partial PINN(x,y,t)}{\partial t}+  \frac{\partial T(y)}{\partial y}\frac{\partial PINN(x,y,t)}{\partial y} - \right . \nonumber \\ \left. 0.1\frac{\partial  PINN(x,y,t)}{\partial x^2}-0.01\frac{\partial  PINN(x,y,t)}{\partial y^2}
 \right)^2
\end{eqnarray}

\begin{eqnarray}
L_{residual}(x,y,t) =
\left(
\frac{\partial u}{\partial t}+  \frac{\partial T}{\partial y}\frac{\partial u}{\partial y} -  K_x\frac{\partial  u}{\partial x^2}- K_y\frac{\partial  u}{\partial y^2} - f
 \right)^2
\end{eqnarray}

We also define the loss for training the initial condition as the residual of the initial condition:

\begin{eqnarray}
L_{Initial}(x,y,0) = 
\left(PINN(x,y,0)-u_0(x,y)\right)^2
\end{eqnarray}

as well as the loss of the residual of the boundary condition:

\begin{eqnarray}
L_{boundary}(x,y,t) = 
f\left(\frac{\partial PINN(x,y,t)}{\partial n}-0\right)^2
\end{eqnarray}

\begin{figure}[!ht]
    \centering
    \includegraphics[width=0.48\textwidth]{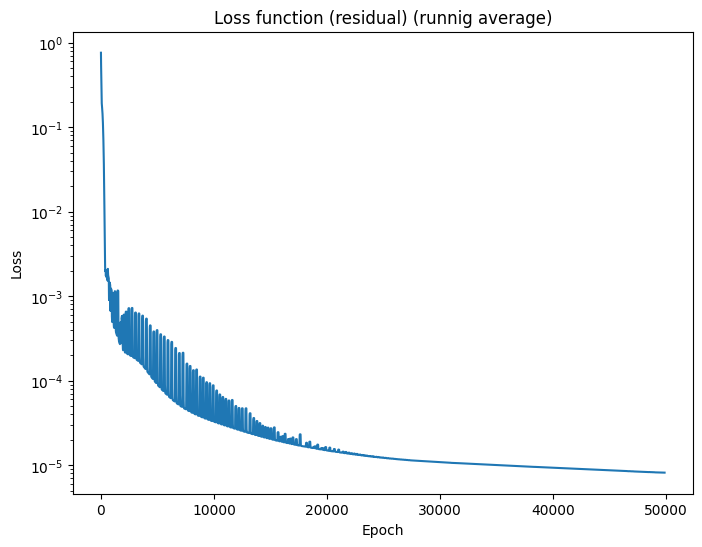}
    \includegraphics[width=0.48\textwidth]{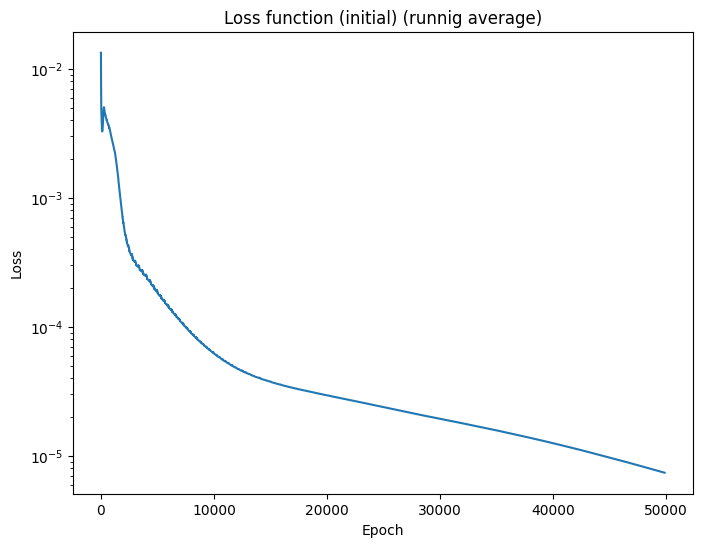}\\
    \includegraphics[width=0.48\textwidth]{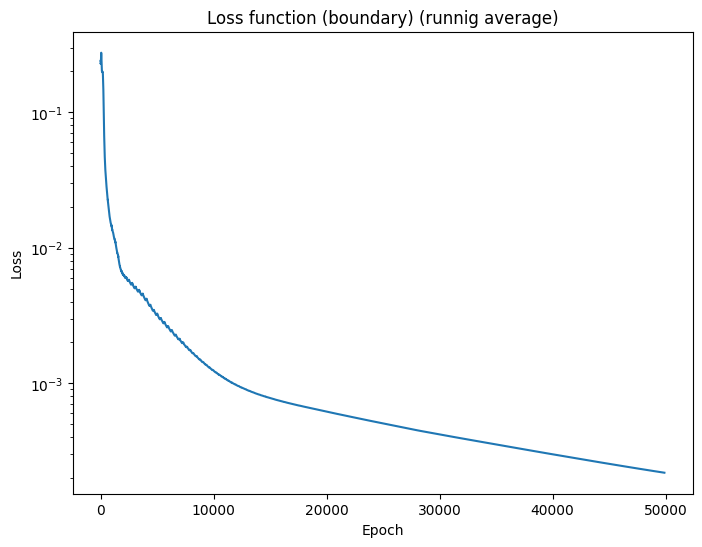}
    \includegraphics[width=0.48\textwidth]{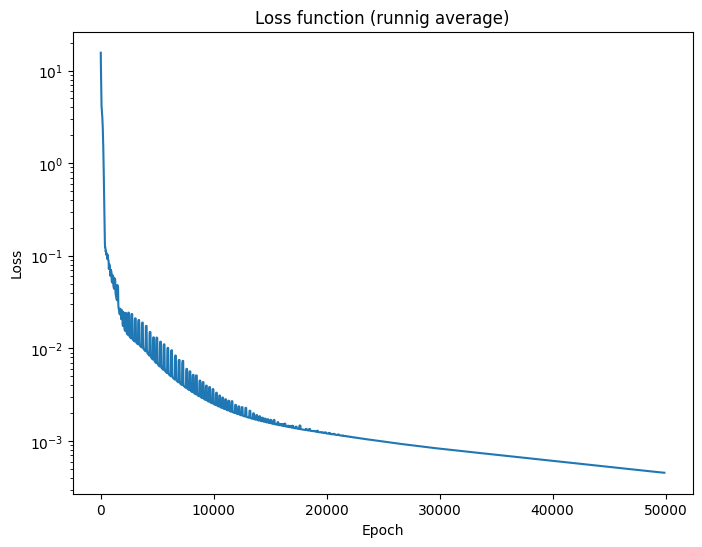}
    \caption{Svalbard winter. The temperature increases in vertical direction close to the ground. The convergence of residual, initial, boundary, and total loss functions.}
    \label{fig:thermal:losses_winter}
\end{figure}

Backpropagation is the core of neural network training, and it employs  the chain rule:

\[ \frac{\partial L}{\partial \mathbf{W}_i} = \frac{\partial L}{\partial \mathbf{a}_i} \cdot \frac{\partial \mathbf{a}_i}{\partial \mathbf{z}_i} \cdot \frac{\partial \mathbf{z}_i}{\partial \mathbf{W}_i} \]

where \( L \) is the loss, \( \mathbf{a}_i \) is the activation, and \( \mathbf{z}_i \) is the input to the activation function at layer \( i \).

The sketch of the training procedure is the following.

\begin{itemize}
    \item \textbf{Gradient Descent (GD)} - is the most basic algorithm, which iteratively adjusts the parameters in the direction of the negative gradient of the loss function to minimize the loss. The update rule for GD is:

    \[ \mathbf{W} \leftarrow \mathbf{W} - \eta \cdot \nabla_{\mathbf{W}} L \]
    \[ \mathbf{b} \leftarrow \mathbf{b} - \eta \cdot \nabla_{\mathbf{b}} L \]
    
    where \( \eta \) is the learning rate, and \( \nabla_{\mathbf{W}} L \) and \( \nabla_{\mathbf{b}} L \) are the gradients of the loss with respect to the weights and biases, respectively.
    
    \item \textbf{Adam optimizer (Adaptive Moment Estimation)} \cite{c34} - is a more advanced optimization algorithm that maintains running averages of both the gradients (first moment) and the squared gradients (second moment). The update rules for Adam are:

    \begin{enumerate}
        \item Compute the exponentially decaying average of past gradients (first moment estimate):
        \[ m_t = \beta_1 m_{t-1} + (1 - \beta_1) g_t \]
    
        \item Compute the exponentially decaying average of past squared gradients (second moment estimate):
        \[ v_t = \beta_2 v_{t-1} + (1 - \beta_2) g_t^2 \]
    
        \item Compute bias-corrected estimates:
        \[ \hat{m}_t = \frac{m_t}{1 - \beta_1^t} \]
        \[ \hat{v}_t = \frac{v_t}{1 - \beta_2^t} \]
    
        \item Update parameters:
        \[ \mathbf{W} \leftarrow \mathbf{W} - \eta \frac{\hat{m}_t}{\sqrt{\hat{v}_t} + \epsilon} \]
        \[ \mathbf{b} \leftarrow \mathbf{b} - \eta \frac{\hat{m}_t}{\sqrt{\hat{v}_t} + \epsilon} \]
    \end{enumerate}

where \( \beta_1 \) and \( \beta_2 \) are hyperparameters that control the decay rates of these running averages, \( \epsilon \) is a small constant to prevent division by zero, and \( g_t \) is the gradient from iteration  \( t \). The general idea of the Adam algorithm is to average the gradients from several past iterations, converging towards global minima and avoiding local minima.
\end{itemize}

\section{The structure of the code}
\label{sec:code}

\subsection{Colab implementation}
\label{sec:colab}

The simulation code may be downloaded from~\url{https://github.com/pmaczuga/pinn-notebooks} and executed in Google Colab in the fully automatic mode.

\subsection{Parameters}
\label{sec:params}

There are the following model parameters that the user can define:
\begin{itemize}
\item {\tt LENGTH}, {\tt TOTAL\_TIME}. The code works in the space-time domain, where the training is performed by selecting point along $x$, $y$
and $t$ axes. The {\tt LENGTH} parameter defines the dimension of the domain along $x$ and $y$ axes. The domain dimension is {\tt [0,LENGTH]x[0,LENGTH]x[0,TOTAL\_TIME]}.
The {\tt TOTAL\_TIME} parameter defines the length of the space-time domain along the $t$ axis. It is the total time of the transient phenomena we want to simulate.
\item {\tt N\_POINTS}. This parameter defines the number of points used for training. By default, the points are selected randomly along $x$, $y$, and $t$ axes. It is easily possible to extend the code to support different numbers of points or different distributions of points along different axes of the coordinate system.
\item {\tt N\_POINTS\_PLOT}. This parameter defines the number of points used to probe the solution and plot the output plots after the training.
\item {\tt WEIGHT\_RESIDUAL}, {\tt WEIGHT\_INITIAL}, {\tt WEIGHT\_BOUNDARY}. These parameters define the weights for the training of residual, initial condition, and boundary condition loss functions.
\item {\tt LAYERS}, {\tt NEURONS\_PER\_LAYER}. These parameters define the neural network by providing the number of layers and number of neurons per neural network layer.
\item {\tt EPOCHS}, and {\tt LEARNING\_RATE} provide a number of epochs and the training rate for the training procedure.
\end{itemize}

During the training, we used the following global parameter values:

\begin{lstlisting}
## Parameters
LENGTH = 1.
TOTAL_TIME = 1.
N_POINTS = 15
N_POINTS_PLOT = 150
WEIGHT_RESIDUAL = 20.0
WEIGHT_INITIAL = 1.0
WEIGHT_BOUNDARY = 10.0
LAYERS = 2
NEURONS_PER_LAYER = 600
EPOCHS = 30_000
LEARNING_RATE = 0.002
\end{lstlisting}

\subsection{PINN class}
\label{sec:pinclass}

The PINN class defines the functionality for a simple neural network accepting three features as input: the values of $(x,y,t)$ and returning a single output, namely the value of the solution $u(x,y,t)$.
We provide the following features:
\begin{itemize}
\item The {\tt f} routine compute the values of the approximate solution at point $(x,y,t)$.
\item 
The routines {\tt dfdt}, {\tt dfdx}, {\tt dfdy} compute the derivatives of the approximate solution at point $(x,y,t)$ with respect to either $x$, $y$, or $t$ using the PyTorch autograd method.
\end{itemize}

We add the definitions of the {\tt Kx} and {\tt Ky} variables into the {\tt Loss} class. 

\subsection{Processing initial and boundary conditions}
\label{sec:bounds}

Since the training is performed in the space-time domain 
{\tt [0,LENGTH]x[0,LENGTH]x} {\tt[0,TOTAL\_TIME]}, we provide in 
\begin{itemize}
    \item 
{\tt get\_interior\_points} the functionality to identify the points from the training of the residual loss, in
\item {\tt get\_initial\_points} the functionality to identify points for the training of the initial loss, and in 
\item {\tt get\_boundary\_points} the functionality for training the boundary loss.
\end{itemize}

\subsection{Loss functions}
\label{sec:loss}

We provide interfaces for defining the loss functions inside the {\tt Loss} class.
Namely, we define the {\tt residual\_loss}, {\tt initial\_loss} and {\tt boundary\_loss}. Since the initial and boundary loss is universal, and residual loss is problem specific, we provide fixed implementations for the initial and boundary losses, assuming that the initial state is prescribed in the {\tt initial\_condition} routine and that the boundary conditions are zero Neumann. The code can be easily extended to support different boundary conditions.

\begin{lstlisting}
class Loss:
...
    def residual_loss(self, pinn: PINN):
        x, y, t = get_interior_points(self.x_domain, self.y_domain, \
            self.t_domain, self.n_points, pinn.device())
        loss = dfdt(pinn, x, y, t).to(device) 
        - self.dTy(y, t)*dfdy(pinn, x, y, t).to(device) 
        - self.Kx*dfdx(pinn, x, y, t,order=2).to(device)  
        - self.Ky*dfdy(pinn, x, y, t, order=2).to(device) 
        - self.source(y,t).to(device)
        return loss.pow(2).mean

    def initial_loss(self, pinn: PINN):
        x, y, t = get_initial_points(self.x_domain, self.y_domain, \ 
            self.t_domain,self.n_points, pinn.device())
        pinn_init = self.initial_condition(x, y)
        loss = f(pinn, x, y, t) - pinn_init
        return loss.pow(2).mean()

    def boundary_loss(self, pinn: PINN):
        down, up, left, right = get_boundary_points(self.x_domain, \ 
            self.y_domain,self.t_domain,self.n_points, pinn.device())
        x_down,  y_down,  t_down    = down
        x_up,    y_up,    t_up      = up
        x_left,  y_left,  t_left    = left
        x_right, y_right, t_right   = right

        L_down  = dfdy( pinn, x_down,  y_down,  t_down  )
        L_up    = dfdy( pinn, x_up,    y_up,    t_up    )
        L_left  = dfdx( pinn, x_left,  y_left,  t_left  )
        L_right = dfdx( pinn, x_right, y_right, t_right )

        return L_down.pow(2).mean()  + \
            L_up.pow(2).mean()    + \
            L_left.pow(2).mean()  + \
            L_right.pow(2).mean()
\end{lstlisting}

The initial condition is defined in the {\tt initial\_condition} routine, which returns a value of the initial condition at point $(x,y,0)$.

\begin{lstlisting}
# Initial condition
def initial_condition(x: torch.Tensor, y: torch.Tensor) -> torch.Tensor:
...
    res = INITIAL POLLUTION DISTRIBUTION AS OBTAINED FROM FEM SOLVER
    return res
\end{lstlisting}

The minimization of the three losses, is the multi-objective optimization problem. The loss functions can be weighted $L=W_{residual}L_{residual}+W_{initial}L_{initial}+W_{boiundary}L_{boundary}$ with the weights $(W_{residual},W_{initial},W_{boiundary})$ selected automatically using the SoftAdapt algorithm \cite{c35}.

The number of neurons and the number of layers in the PINNs can be estimated using the results of Jinchao Xu, showing the  analogies between neural networks and linear and higher-order finite element methods \cite{c36,c37}. The weights of the loss functions for the multi-objective optimization can be determined automatically using SoftAdapt algorithm \cite{c35}.

\subsection{Summer simulation}
\label{sec:summer}

In this section, we present numerical results of the pollution dissipation computed for the vertical temperature profile during the summer day.

In summer, temperature inversions are less common on the Spitzbergen than in winter. The temperature in the troposphere (lower layer of the atmosphere) usually decreases with altitude.

\begin{figure}[!ht]
      \centering
      \includegraphics[width=0.32\textwidth]{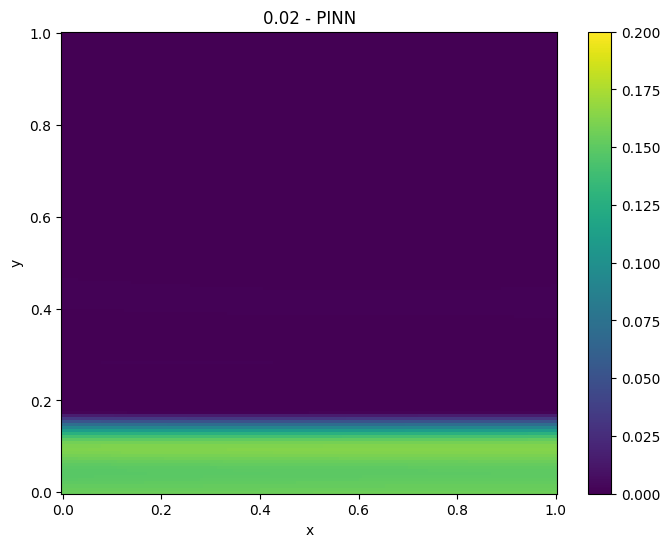}
      \includegraphics[width=0.32\textwidth]{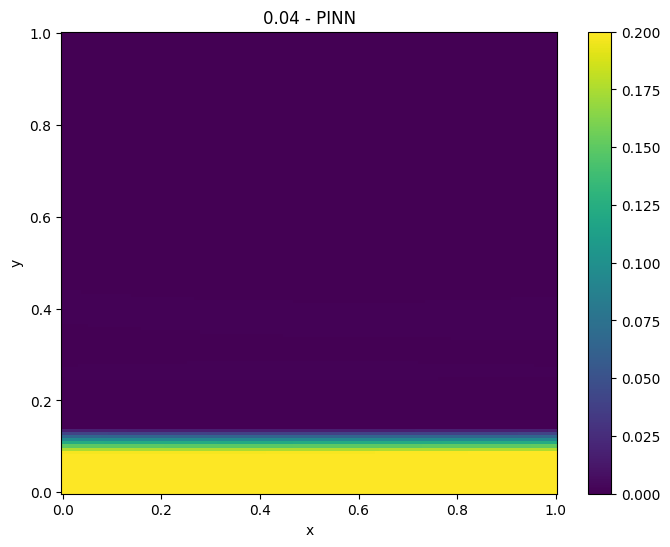}
      \includegraphics[width=0.32\textwidth]{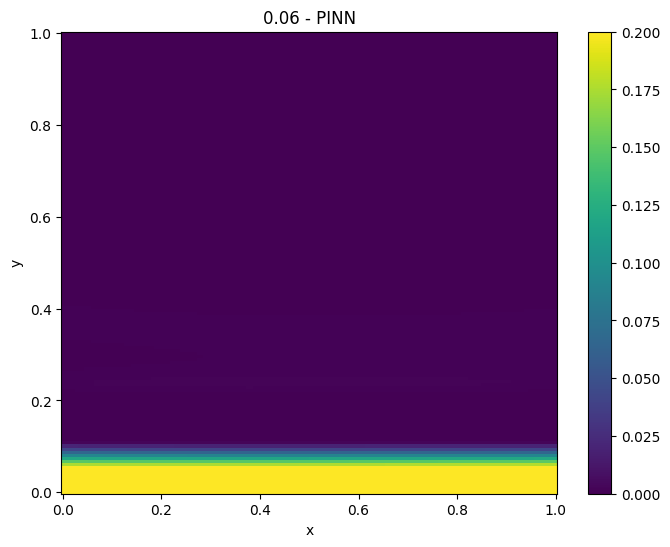}\\
       \includegraphics[width=0.32\textwidth]{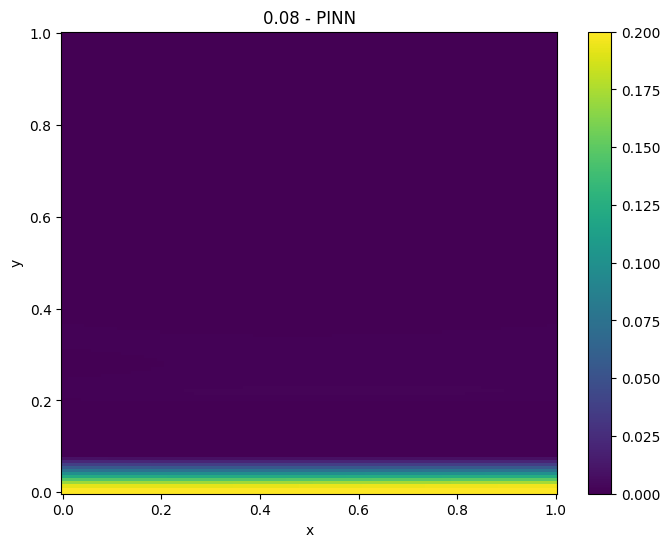}
      \includegraphics[width=0.32\textwidth]{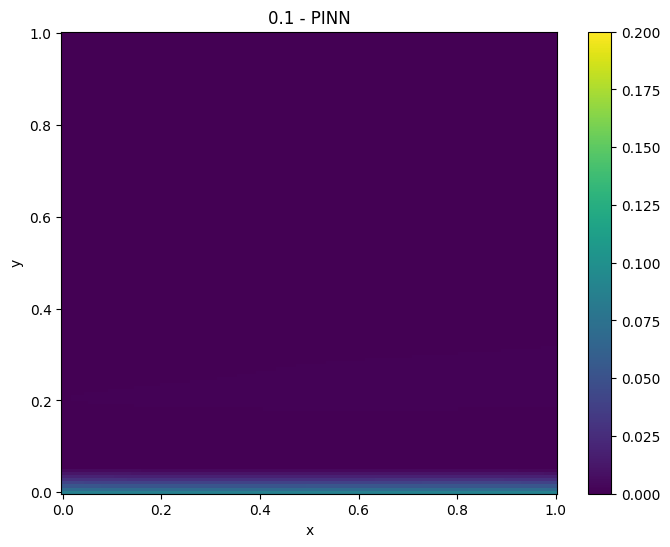}
      \includegraphics[width=0.32\textwidth]{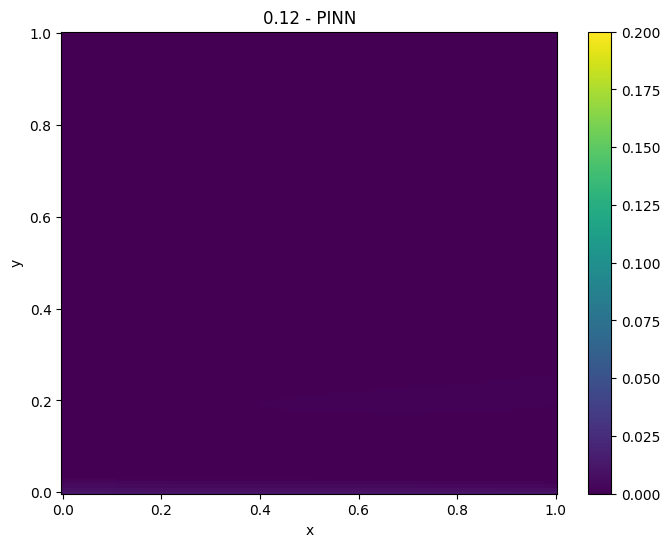}\\
      \caption{Pollution concentration during the Svalbard summer where the temperature decreases in the vertical direction.}
      \label{fig:thermal_summer}
\end{figure}

The convergence of the loss functions is summarized in Figure~\ref{fig:thermal:losses_summer}. 
The snapshots from the simulations are presented in Figure~\ref{fig:thermal_summer}.
The pollution concentration units are dimensionless, and the goal of the simulation is to present the quantitative behavior of the pollution propagation with the temperature profile during the summer period.
The pollution generated by the electric power plant 
dissipates due to the vertical temperature gradients.

\begin{figure}
    \centering
    \includegraphics[width=0.3\textwidth]{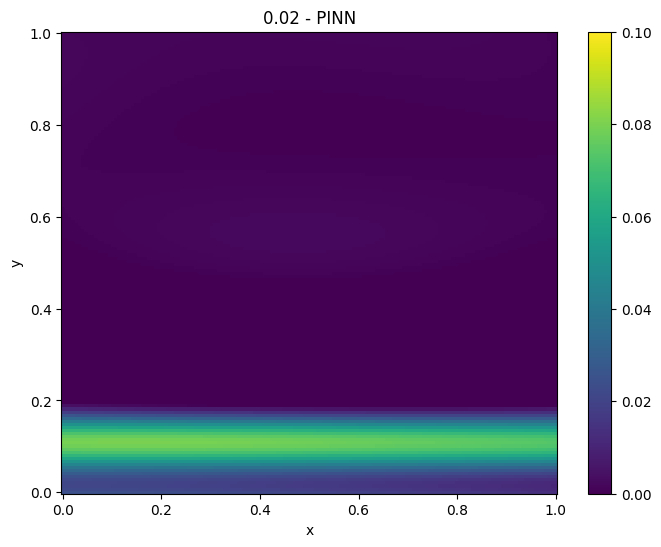}
    \includegraphics[width=0.3\textwidth]{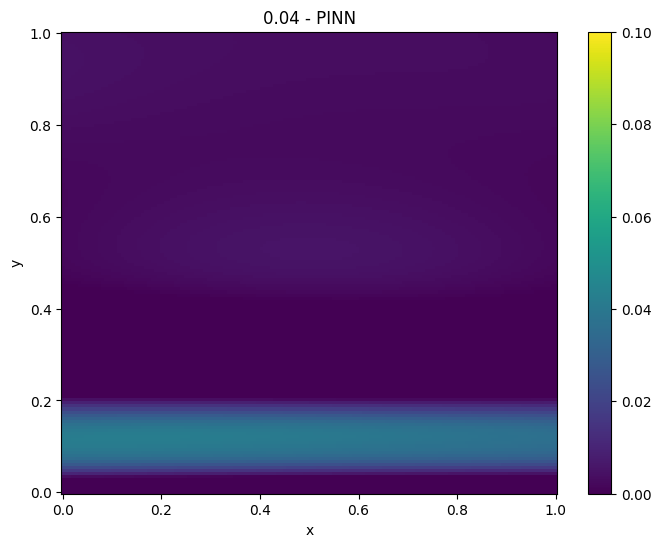}
    \includegraphics[width=0.3\textwidth]{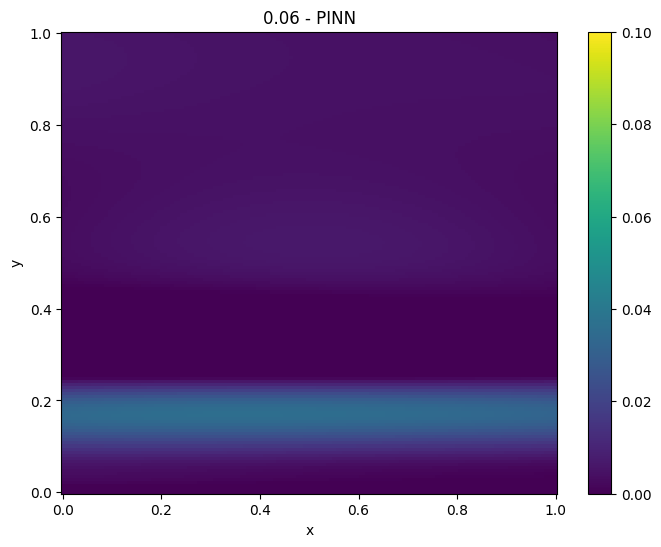}\\ 
   \includegraphics[width=0.3\textwidth]{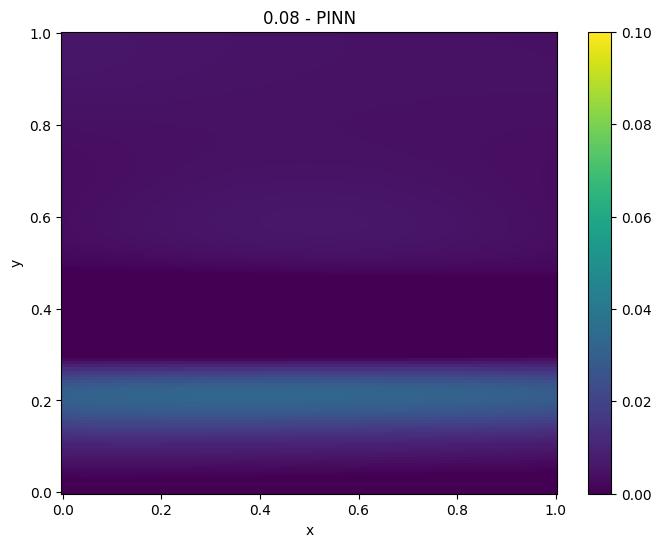}
    \includegraphics[width=0.3\textwidth]{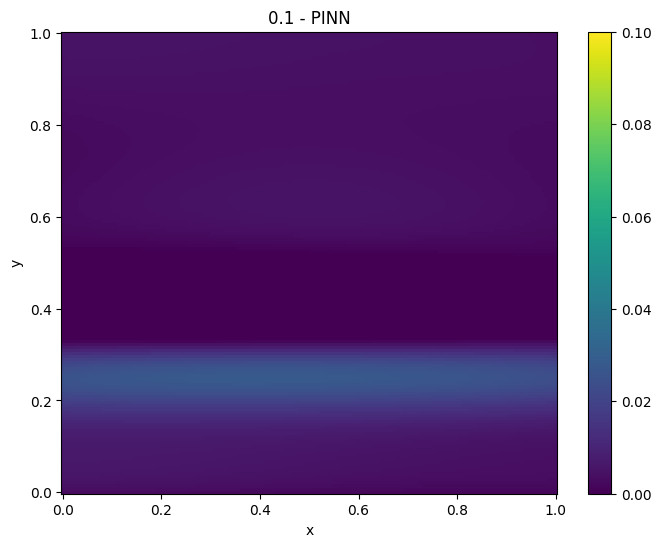}
    \includegraphics[width=0.3\textwidth]{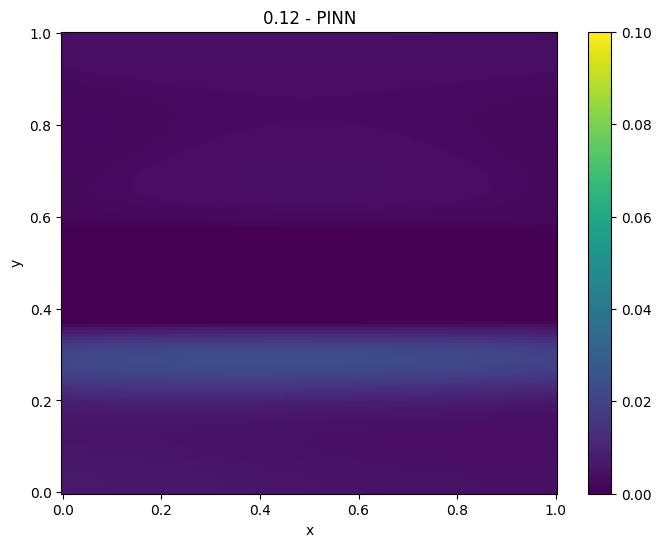}\\
    \caption{Thermal inversion simulation for the Svalbard winter, where the temperature increases in the vertical direction.}
    \label{fig:thermal_winter}
\end{figure}

\subsection{Winter simulation}
\label{sec:winter}

In this section, we present numerical results of the pollution dissipation computed for the temperature profiles during the winter night.

The convergence of the loss functions is summarized in Figure~\ref{fig:thermal:losses_winter}. 
The snapshots from the simulations are presented in Figure~\ref{fig:thermal_winter}.

The dimensionless pollution concentration units illustrate the quantitative behavior of the pollution propagation with the temperature profile from the Arctic night.

Just like in urban environments, thermal inversions in the Arctic can trap pollutants.
The absence of sunlight during the Arctic night leads to intense cooling of the Earth's surface. As the ground loses heat, the air directly above it also cools rapidly.
The lack of solar heating during the Arctic night results in stable atmospheric conditions with minimal vertical air mixing. This stability allows the cold air to remain trapped near the surface.

\subsection{Notes on computational cost}
The two-dimensional time-dependent PINN simulator execution on A100 Backend Google Compute Engine with Python 3 and GPU graphic card equipped with 83.48 GB of memory and 235.68 GB disc space takes
around 15 minutes of computing time.
This execution time is comparable with the execution of 200 times steps of the non-stationary three-dimensional graph-grammar-based finite element method solver
on a laptop with with 11th Gen Intel(R) Core(TM) i5-11500H @ 2.90GHz, 2.92 GHz, and 32 GB of RAM,
providing an estimate of 9 hours of real-time pollution generation from a chimney.
The finite element method simulation requires the development of the stabilized time-integration scheme, in our case, the Crank-Nicolson method.
The PINNs do not require the development of a stabilized time-integration scheme; the time-dependent problem is trained in the space-time domain. 

\section{Conclusions}
\label{sec:conclus}

We presented an original model describing the production of graph grammars, transformation sequences of graphs representing a computational grid, expressing an algorithm for adapting a three-dimensional computational mesh that does not generate hanging nodes.
The graph transformation rules model the Rivara algorithms. Its
the idea is to break elements along the longest edges and propagate the refinement to adjacent elements to avoid hanging nodes in three-dimensional computational grids. 
The graph transformations were used to generate a computational grid for simulating pollution propagation from a coal-fired power plant in Longyearbyen, Spitzbergen.
We also introduce a computational code performing Physics Informed
Neural Networks simulations of the pollution propagation. The PINNs are attractive alternatives for simulations carried out using the finite element method. 
They do not require a time integration scheme and do not generate stability problems encountered by time integration.
On the other hand, successful training of the PINN model is a multi-objective optimization problem, and it requires guessing several model parameters, such as the number of layers of the neural network, the number of neurons, the training rate, and the loss function weights for training.
Some analogies between neural networks and linear and higher-order finite element methods can be found in works of \cite{c36,c37}. They enable us to estimate the size of the neural networks. The weights of the loss functions for the multi-objective optimization can be determined automatically using the SoftAdapt algorithm \cite{c35}. Nevertheless, the actual state-of-the-art PINNs enable, in the authors' opinion, the successful and efficient application of the two-dimensional PINN model in engineering applications.

\section{Acknowledgments}
The Authors gratefully acknowledge the support and assistance of The Polish Polar Station Hornsund for help with data collection.

The authors are grateful for support from the funds the Polish Ministry of Science and Higher Education assigned to AGH University of Krakow.
The work supported by  ``Excellence initiative - research university" for the AGH University of Krakow.

The work of Albert Oliver-Serra was supported by 
"Ayudas para la recualificación del sistema universitario español" grant funded by the ULPGC, the Ministry of Universities by Order UNI/501/2021 of 26 May, and the European Union-Next Generation EU Funds.


 \bibliographystyle{elsarticle-num} 
 \bibliography{cas-refs}

\begin{thebibliography}{00}

\bibitem{c1}
K. Podsiadło, A. O. Serra, A. Paszyńska, R. Montenegro, I. Henriksen,
M. Paszyński, K. Pingali, Parallel graph-grammar-based algorithm for the
longest-edge refinement of triangular meshes and the pollution simulations in
lesser poland area, Engineering with Computers 37 (2021) 3857–3880.
\bibitem{c2}
A. N. Brooks, T. J. Hughes, Streamline upwind/petrov-galerkin formulations
for convection dominated flows with particular emphasis on the incompressible
navier-stokes equations, Computer Methods in Applied Mechanics and Engineering
32 (1) (1982) 199–259. doi:https://doi.org/10.1016/0045-7825(82)90071-8.
URL https://www.sciencedirect.com/science/article/pii/
0045782582900718
\bibitem{c3} M. C. Rivara, Algorithms for refining triangular grids suitable for adaptive and
multigrid techniques, International Journal Numerical Methods in Engineering
20 (4) (1984) 745–756.
\bibitem{c4}
M. C. Rivara, Mesh refinement processes based on the generalized bisection of
simplices, SIAM Journal Numerical Analysis 21 (3) (1984) 604–613.
\bibitem{c5}
A. Paszyńska, M. Paszyński, E. Grabska, Graph transformations for modeling
hp-adaptive finite element method with triangular elements, in: Proceedings of
the 8th International Conference on Computational Science, Part III, ICCS ’08,
Springer-Verlag, Berlin, Heidelberg, 2008, p. 604–613. doi:10.1007/978-3-540-
69389-5\_68.
URL https://doi.org/10.1007/978-3-540-69389-5\_68
\bibitem{c6}
 A. Paszyńska, M. Paszyński, E. Grabska, Graph transformations for modeling
hp-adaptive finite element method with mixed triangular and rectangular
elements, in: Proceedings of the 9th International Conference on Computational
Science, ICCS 2009, Springer-Verlag, Berlin, Heidelberg, 2009, p. 875–884.
doi:10.1007/978-3-642-01973-9\_97.
URL https://doi.org/10.1007/978-3-642-01973-9\_97
\bibitem{c7}
M. Paszyński, On the parallelization of self-adaptive hp-finite element methods
part i. composite programmable graph grammar model, Fundamenta Informaticae
93 (4) (2009) 411–434.
\bibitem{c8}
 M. Paszyński, On the parallelization of self-adaptive hp-finite element methods
part ii. partitioning communication agglomeration mapping (pcam) analysis,
Fundamenta Informaticae 93 (4) (2009) 435–457.

\bibitem{c9}
W. B. F. Ryan, S. M. Carbotte, J. O. Coplan, S. O’Hara, A. Melkonian,
R. Arko, R. A. Weissel, V. Ferrini, A. Goodwillie, F. Nitsche,
J. Bonczkowski, R. Zemsky, Global multi-resolution topography
synthesis, Geochemistry, Geophysics, Geosystems 10 (3) (2009).
arXiv:https://agupubs.onlinelibrary.wiley.com/doi/pdf/10.1029/2008GC002332,
doi:https://doi.org/10.1029/2008GC002332.
URL https://agupubs.onlinelibrary.wiley.com/doi/abs/10.1029/
2008GC002332
\bibitem{c10}
M. Raissi, P. Perdikaris, G. Karniadakis, Physics-informed neural networks: A
deep learning framework for solving forward and inverse problems involving
nonlinear partial differential equations, Journal of Computational Physics 378
(2019) 686–707. doi:https://doi.org/10.1016/j.jcp.2018.10.045.
\bibitem{c11}
G. Hinton, L. Deng, D. Yu, G. E. Dahl, A.-r. Mohamed, N. Jaitly, A. Senior,
V. Vanhoucke, P. Nguyen, T. N. Sainath, et al., Deep neural networks for acoustic
modeling in speech recognition: The shared views of four research groups, IEEE
Signal processing magazine 29 (6) (2012) 82–97.
\bibitem{c12}
A. Krizhevsky, I. Sutskever, G. E. Hinton, Imagenet classification with deep
convolutional neural networks, Communications of the ACM 60 (6) (2017) 84–90.
\bibitem{c13}
M. Gheisari, G. Wang, M. Z. A. Bhuiyan, A survey on deep learning in big data,
in: 2017 IEEE international conference on computational science and engineering
(CSE) and IEEE international conference on embedded and ubiquitous computing
(EUC), Vol. 2, IEEE, 2017, pp. 173–180.
\bibitem{c14}
M. Raissi, P. Perdikaris, G. E. Karniadakis, Physics-informed neural networks:
A deep learning framework for solving forward and inverse problems involving
nonlinear partial differential equations, Journal of Computational physics 378
(2019) 686–707.
\bibitem{c15}
S. Cai, Z. Mao, Z. Wang, M. Yin, G. E. Karniadakis, Physics-informed neural
networks (PINNs) for fluid mechanics: A review, Acta Mechanica Sinica 37 (12)
(2021) 1727–1738.
\bibitem{c16}
Z. Mao, A. D. Jagtap, G. E. Karniadakis, Physics-informed neural networks for
high-speed flows, Computer Methods in Applied Mechanics and Engineering 360
(2020) 112789.

\bibitem{c17}
J. Ling, A. Kurzawski, J. Templeton, Reynolds averaged turbulence modelling
using deep neural networks with embedded invariance, Journal of Fuild Mechanics
807 (2016) 155–166. doi:10.1017/jfm.2016.615.
\bibitem{c18}
L. Sun, H. Gao, S. Pan, J.-X. Wang, Surrogate modeling for fluid flows based on
physics-constrained deep learning without simulation data, Computer Methods in
Applied Mechanics and Engineering 361 (2020). doi:10.1016/j.cma.2019.112732.
\bibitem{c19}
N. Wandel, M. Weinmann, M. Neidlin, R. Klein, Spline-pinn: Approaching
pdes without data using fast, physics-informed hermite-spline cnns, Proceedings
of the AAAI Conference on Artificial Intelligence 36 (8) (2022) 8529–8538.
doi:10.1609/aaai.v36i8.20830.
URL https://ojs.aaai.org/index.php/AAAI/article/view/20830
\bibitem{c20}
M. Rasht-Behesht, C. Huber, K. Shukla, G. E. Karniadakis, Physics-informed
neural networks (pinns) for wave propagation and full waveform inversions,
Journal of Geophysical Research: Solid Earth 127 (5) (2022) e2021JB023120.
\bibitem{c21}
N. Geneva, N. Zabaras, Modeling the dynamics of pde systems with physicsconstrained
deep auto-regressive networks, Journal of Computational Physics
403 (2020). doi:10.1016/j.jcp.2019.109056.
\bibitem{c22}
S. Goswami, C. Anitescu, S. Chakraborty, T. Rabczuk, Transfer learning
enhanced physics informed neural network for phase-field modeling
of fracture, Theoretical and applied fracture machanics 106 (2020).
doi:10.1016/j.tafmec.2019.102447.
\bibitem{c23}
M. Alber, A. B. Tepole, W. R. Cannon, S. De, S. Dura-Bernal, K. Garikipati,
G. Karniadakis, W. W. Lytton, P. Perdikaris, L. Petzold, E. Kuhl, Integrating
machine learning and multiscale modeling-perspectives, challenges, and opportunities
in the biologica biomedical, and behavioral sciences, NPJ Digital Medicine
2 (2019). doi:10.1038/s41746-019-0193-y.
\bibitem{c24}
G. Kissas, Y. Yang, E. Hwuang, W. R. Witschey, J. A. Detre, P. Perdikaris,
Machine learning in cardiovascular flows modeling: Predicting arterial blood
pressure from non-invasive 4d flow mri data using physics-informed neural
networks, Computer Methods in Applied Mechanics and Engineering 358 (2020).
doi:10.1016/j.cma.2019.112623.
\bibitem{c25}
H. Jin, M. Mattheakis, P. Protopapas, Physics-informed neural networks for
quantum eigenvalue problems, in: 2022 International Joint Conference on Neural
Networks (IJCNN), 2022, pp. 1–8. doi:10.1109/IJCNN55064.2022.9891944.

\bibitem{c26}
R. Nellikkath, S. Chatzivasileiadis, Physics-informed neural networks for
minimising worst-case violations in dc optimal power flow, in: 2021
IEEE International Conference on Communications, Control, and Computing
Technologies for Smart Grids (SmartGridComm), 2021, pp. 419–424.
doi:10.1109/SmartGridComm51999.2021.9632308.
\bibitem{c27}
X. Huang, H. Liu, B. Shi, Z. Wang, K. Yang, Y. Li, M. Wang, H. Chu, J. Zhou,
F. Yu, B. Hua, B. Dong, L. Chen, A universal pinns method for solving partial differential
equations with a point source, Proceedings of the Fourteen International
Joint Conference on Artificial Intelligence (IJCAI-22) (2022) 3839–3846.
\bibitem{c28}
Y. Yang, P. Perdikaris, Adversarial uncertainty quantification in physicsinformed
neural networks, Journal of Computational Physics 394 (2019) 136–152.
doi:10.1016/j.jcp.2019.05.027.
\bibitem{c29}
F. Sun, Y. Liu, H. Sun, Physics-informed spline learning for nonlinear dynamics
discovery, Proceedings of the Thirtieth International Joint Conference on
Artificial Intelligence (IJCAI-21) (2021) 2054–2061.
\bibitem{c30}
J. Kim, K. Lee, D. Lee, S. Y. Jin, N. Park, Dpm: A novel training method for
physics-informed neural networks in extrapolation, in: 35th AAAI Conference
on Artificial Intelligence, AAAI 2021, 2021, pp. 8146–8154.
\bibitem{c31}
Y. Chen, L. Lu, G. E. Karniadakis, L. Dal Negro, Physics-informed neural
networks for inverse problems in nano-optics and metamaterials, Optics express
28 (8) (2020) 11618–11633.
\bibitem{c32}
S. Mishra, R. Molinaro, Estimates on the generalization error of physics-informed
neural networks for approximating a class of inverse problems for PDEs, IMA
Journal of Numerical Analysis 42 (2) (2022) 981–1022.
\bibitem{c33}
L. Lu, R. Pestourie, W. Yao, Z. Wang, F. Verdugo, S. G. Johnson, Physicsinformed
neural networks with hard constraints for inverse design, SIAM Journal
on Scientific Computing 43 (6) (2021) B1105–B1132. doi:10.1137/21M1397908.
\bibitem{c34}
D. P. Kingma, J. Ba, Adam: A method for stochastic optimization, arXiv
preprint arXiv:1412.6980 (2014).
\bibitem{c35}
A. A. Heydari, C. A. Thompson, A. Mehmood, Softadapt: Techniques for
adaptive loss weighting of neural networks with multi-part loss functions (2019).
arXiv:1912.12355.
URL https://arxiv.org/abs/1912.12355

\bibitem{c36}
J. He, L. Li, J. Xu, C. Zheng, J. He, Relu deep neural networks and linear
finite elements, Journal of Computational Mathematics 38 (3) (2020) 502–527.
doi:10.4208/jcm.1901-m2018-0160.
URL http://dx.doi.org/10.4208/jcm.1901-m2018-0160
\bibitem{c37}
J. He, J. Xu, Deep neural networks and finite elements of any order on arbitrary
dimensions (2024). arXiv:2312.14276.
URL https://arxiv.org/abs/2312.14276

 \end{thebibliography}


\end{document}